\documentclass[aap,MSNbibl,dvips]{arximspdf}
\usepackage{mathrsfs}
\usepackage{graphicx}
\doi{10.1214/13-AAP935} 
\volume{24}
\issue{2}
\pubyear{2014}
\firstpage{760}
\lastpage{810}

\makeatletter

\newcommand{\rrvert}{\vert}
\newcommand{\llvert}{\vert}
\def\cal{\mathcal}

\newcommand{\argmin}{\operatorname{arg\,min}}
\newcommand{\argmax}{\operatorname{arg\,max}}

\newtheorem{lemma}{Lemma}[section]
\newtheorem{theorem}{Theorem}[section]

\newproclaim{definition}{Definition}[section]
\newproclaim{example}{Example}[section]
\newtheorem{proposition}{Proposition}[section]

\newproclaim{assumption}{Assumption}[section]

\newproclaim{remark}{Remark}[section]

\newcommand{\eps}{\varepsilon}
\newcommand{\Ex}{\mathbb{E}}

\newcommand{\vr}{\varrho}
\newcommand{\al}{\alpha}
\newcommand{\s}{\sigma}
\newcommand{\mD}{\mathcal{D}}
\newcommand{\mB}{\mathcal{B}}

\newcommand{\Gam}{{\Gamma}}
\newcommand{\Del}{{\Delta}}
\newcommand{\La}{{\Lambda}}
\newcommand{\X}{{\Xi}}
\newcommand{\PI}{{\Pi}}
\newcommand{\Sig}{{\Sigma}}
\newcommand{\Om}{{\Omega}}
\newcommand{\D}{\mathbb{D}}
\newcommand{\N}{\mathbb{N}}
\newcommand{\R}{\mathbb{R}}
\newcommand{\hX}{\hat{X}}
\newcommand{\hB}{\hat{B}}
\newcommand{\hI}{\hat{I}}

\newcommand{\hA}{\hat{A}}
\newcommand{\tB}{\tilde{B}}

\newcommand{\cX}{\check{X}}
\newcommand{\cQ}{\check{Q}}

\newcommand{\hQ}{\hat{Q}}
\newcommand{\EE}{\mathbb{E}}
\newcommand{\PP}{\mathbb{P}}

\newcommand{\calE}{{\cal E}}
\newcommand{\calF}{{\cal F}}
\newcommand{\calG}{{\cal G}}
\newcommand{\calI}{{\cal I}}
\newcommand{\calJ}{{\cal J}}
\newcommand{\calL}{{\cal L}}
\newcommand{\calP}{{\cal P}}
\newcommand{\calT}{{\cal T}}
\newcommand{\calV}{{\cal V}}

\newcommand{\bI}{{\mathbf I}}
\newcommand{\bJ}{{\mathbf J}}
\newcommand{\bK}{{\mathbf K}}

\newcommand{\scrM}{\mathscr{M}}

\newcommand{\eqref}[1]{(\ref{#1})}

\newcommand{\w}{\wedge}
\newcommand{\To}{\Rightarrow}
\newcommand{\til}{\tilde}
\newcommand{\wh}{\hat}
\newcommand{\iy}{\infty}
\newcommand{\IA}{\mathit{IA}}
\newcommand{\ibar}{{\bar{\imath}}}
\newcommand{\jbar}{{\bar{\jmath}}}

\newcommand{\I}{\mathcal{I}}

\newcommand{\J}{\mathcal{J}}
\newcommand{\mL}{\mathcal{L}}

\newcommand{\tinf}{\rightarrow\infty}
\newcommand{\Pd}{\mathbb{P}}

\makeatother

\begin{document}
\begin{frontmatter}

\title{Scheduling parallel servers in the nondegenerate slowdown
diffusion regime: Asymptotic~optimality~results\thanksref{T1}}
\thankstext{T1}{Supported in part by
the US--Israel BSF (Grant 2008466), the ISF (Grant 1349/08) and the Technion's
fund for promotion of research.}
\runtitle{Asymptotically optimal scheduling}

\begin{aug}
\author[A]{\fnms{Rami} \snm{Atar}\corref{}\ead[label=e1]{atar@ee.technion.ac.il}}
\and
\author[B]{\fnms{Itai} \snm{Gurvich}}
\runauthor{R. Atar and I. Gurvich}
\affiliation{Technion--Israel Institute of Technology and Northwestern
University}
\address[A]{Department of Electrical Engineering\\
Technion--Israel Institute of Technology\\
Haifa 32000\\
Israel} 
\address[B]{Kellogg School of Management\\
Northwestern University\\
Evanston, Illinois 60201\\
USA}
\end{aug}

\received{\smonth{11} \syear{2011}}
\revised{\smonth{12} \syear{2012}}

%
\begin{abstract}
We consider the problem of minimizing queue-length costs in a system
with heterogenous parallel servers, operating in a many-server
heavy-traffic regime with nondegenerate slowdown. This regime is
distinct from the well-studied heavy traffic diffusion regimes, namely
the (single server) conventional regime and the (many-server)
Halfin--Whitt regime. It has the distinguishing property that waiting
times and service times are of comparable magnitudes. We establish an
asymptotic lower bound on the cost and devise a sequence of policies
that asymptotically attain this bound. As in the conventional regime,
the asymptotics can be described by means of a Brownian control problem,
the solution of which exhibits a state space collapse.
\end{abstract}

%
\pdfkeywords{60K25, 60J60, 60F17, 90B22, 68M20, The parallel server model,
many-server queues, heavy traffic, diffusion limits, asymptotically optimal control, nondegenerate slowdown regime}

\end{frontmatter}

\section{Introduction} \label{sec1}

Many-server approximations are ubiquitous in the modeling of
large-scale service systems. A prevalent mode of analysis in this
context is the \textit{Halfin--Whitt} heavy traffic diffusion regime
\cite
{halwhi}, also called
the \textit{quality and efficiency driven} (QED) regime \cite{GMR}. For
the $M/M/N$ queue, a sequence of systems in heavy traffic (HT), indexed
by $n$, is constructed by letting the number of servers, $N^n$, and the
arrival rate, $\lambda^n$, grow with $n$ while the service rate $\mu$
remains fixed, so that the utilization in the $n$th system, $\rho
^n:=\lambda
^n/(N^n\mu^n)$, behaves like
\[
\rho^n=1-O \biggl(\frac{1}{\sqrt{\lambda^n}} \biggr)=1-O \biggl(
\frac
{1}{\sqrt{N^n}} \biggr).
\]
Customer waiting times in this regime are of the order $1/\sqrt {\lambda
^n}$ and are thus order of magnitudes smaller than the service times.
It has been argued that this order of magnitude relation renders the
analysis in this regime relevant for some call centers and certain
health-care systems to whose study it has been applied; see, for
example, \cite{armony2010fair,yuliathesis}. The Halfin--Whitt regime
is typically contrasted with the so-called \textit{conventional} HT
diffusion regime.
Conventional limits are obtained by fixing the number of servers
(typically~1) and letting both the arrival and service rate scale so
that
\[
\rho^n=1-O \biggl(\frac{1}{\sqrt{\lambda^n}} \biggr)=1-O \biggl(
\frac
{1}{\sqrt{\mu^n}} \biggr).
\]
In this regime, the waiting time is of the order of $1/\sqrt{\mu^n}$ so
that, in perfect contrast with the Halfin--Whitt regime, the service
time is negligible compared to the waiting time.

From a modeling viewpoint, it is desirable
to allow for these two important performance measures to be comparable
under the scaling.
Gurvich, Mandelbaum and Shaikhet \cite{man,mansha,gur}
and, independently, Whitt \cite{whi2} have identified a many-server
regime with this property.
Limits for the $M/M/N$ queue in this regime are obtained by scaling the
parameters
so that $\mu^n$ and $N^n$ are of the order of $\sqrt{\lambda^n}$ and
%
%
\begin{equation}
\rho^n=1-O \biggl(\frac{1}{\sqrt{\lambda^n}} \biggr)=1-O \biggl(
\frac
{1}{\sqrt{\mu^nN^n}} \biggr) =1-O \biggl(\frac{1}{N^n} \biggr).\label{eq:HTcond}
\end{equation}
One defines the \textit{slowdown} of a queueing system
as the ratio between the sojourn time and the service time
experienced by a typical customer.
By the foregoing discussion among the three regimes alluded to above,
regime \eqref{eq:HTcond}
is unique in that the slowdown does not degenerate, in the sense that
it does not converge to one of the extreme values, 1 or $\iy$.
We therefore refer to it as the \textit{nondegenerate slowdown} (NDS)
diffusion regime.
This term was coined in \cite{RamiNDS}, where a queue with heterogenous
servers was
analyzed in this regime, and the limiting joint law of waiting time and
service time
was identified.
Both the conventional and the NDS regimes are often referred to as \textit{efficiency driven} (ED).
We refer the reader to \cite{RamiNDS} for further discussion of the
three regimes and to
\cite{GMR} for the distinction between QED and ED regimes. The
relevance of the NDS diffusion regime in real-world applications has
been argued
in \cite{RamiNDS} by demonstrating that some call centers do operate
with comparable
delays and service times; particularly, this is the case for
the detailed empirical study of a call center performed in \cite{man-stat}.
This makes a strong case for analyzing these systems by NDS diffusion
approximations
(see \cite{RamiNDS} for further discussion on this modeling issue, and
Whitt \cite{whi3} for an alternative (ED) regime with comparable delays
and service times).

Control of queueing networks under both the conventional and Halfin--Whitt
diffusion regimes (as well as
fluid regimes) has been an active research area in recent years.
Particularly, the \textit{parallel server model} has been studied in
this context,
where customers of a number of classes
are served in parallel by servers of various types, and a system
administrator dynamically controls
the routing. (For the conventional diffusion regime
see \cite{harlop,belwil1,belwil2,mansto} and for the Halfin--Whitt
regime see, e.g., \cite{stolyar2010control} and
references therein.) In this paper we study the problem of minimizing
queue-length costs in a parallel server model with renewal arrival
processes and exponential service times operating in the NDS regime.
From a control standpoint, a distinctive property of the NDS regime is
that \textit{sojourn time} cost criteria are meaningful as neither the
service time or waiting-time degenerate asymptotically. Having solved
the queue-length problem, we argue (heuristically) how in a simple
case, the sojourn time problem can be reduced to a queue-length problem.
This provides further motivation for the latter.
We leave open the rigorization of this argument and the question of how
general this reduction is, as well as the extension of this work to
general service-time distributions.

In terms of the asymptotic behavior, the NDS control problem is close
to the one in the
conventional regime.
In particular, the Brownian control problem (BCP)
which describes the asymptotics is the same as the one
studied in \cite{harlop} for the conventional regime.
This BCP, under a \textit{complete resource pooling} (CRP) condition,
undergoes a reduction to a one-dimensional problem. This reduction is
often called
a \textit{state space collapse}.
Bell and Williams \cite{belwil1,belwil2} and Mandelbaum and
Stolyar \cite{mansto}
have studied the parallel server model in conventional heavy traffic
and obtained asymptotic
optimality results under the CRP condition.
Bell and Williams address linear costs and construct certain threshold
type policies.
Mandelbaum and Stolyar consider separable convex (and flat at the origin)
costs and work with policies that obey the generalized-$c\mu$ rule.

In this paper we aim at a relatively general cost of the form $\int_0^uC(\hat Q^n(t))\,dt$, where we denote by $\hat Q^n$ a properly scaled
version of queue-length, and use the term ``cost'' to mean a random
variable that is to be minimized stochastically. The function $C$,
referred to as the cost function, satisfies an assumption slightly
weaker than convexity (Assumption \ref{asum:Cstar}), as well as an
assumption regarding the existence of a continuous minimizer
(Assumption \ref{asum:continuousselection}). The first main result
(Theorem \ref{th1}) asserts that the BCP value constitutes an asymptotic
lower bound on the cost under any sequence of policies. The second main result
(Theorem \ref{thm:conv_track}) shows that this lower bound is tight.
The price paid for the generality of $C$ is that $c\mu$ type policies
must be abandoned, and
a more complicated policy has to be used to attain the lower bound.
We introduce a \textit{tracking} type policy, in which
a certain target process [denoted by $\check X^n$ in \eqref{86}] is
computed, and
routing is performed so as to keep the difference
between the actual queue-length process $\hat Q^n$ and the target, small.
It turns out that the techniques required here are
quite different from those in the conventional HT,
due both to the general cost and the difference between the regimes.

Theorem \ref{th1} provides a lower bound that is weaker than what, in
the conventional regime,
is referred to as a \textit{pathwise} bound (as, e.g., in \cite{mansto}).
As we discuss in Section~\ref{sec:model}, in the NDS regime a corresponding
pathwise lower bound does not hold,
and this is the main reason for the complicated proof of the result
as compared to that of the pathwise lower bound in the conventional regime.

The main part of the proof of asymptotic optimality\vspace*{1pt} of the proposed
policy is
showing that the difference $\hat Q^n-\check X^n$ is small. (From that,
asymptotic optimality
follows rather easily because the process $\hat Q^n$ can then be shown
to behave like the
explicit solution to the BCP.)
This proof is based on showing that re-balancing of the workload among
the queues can occur quickly
on the relevant time scale. This is facilitated by the fact that
the workload is evenly divided (in the sense of order of magnitude)
between the queues and the servers and
short service times allow to move significant workload from one queue
to the other before
the total-workload process changes considerably.
(This explains the similarity to the behavior in conventional heavy
traffic, where the same quick
response is possible.
On the other hand, such nearly instantaneous re-balancing cannot be
performed in the
Halfin--Whitt regime in which most of the workload is in service and,
to re-direct a nonnegligible
fraction of workload from one class to the other requires a large
number of service completions.
Indeed, the resulting Brownian control problem there is different \cite
{atar2005scheduling}.)

An analysis of the problem in the case of homogenous servers and interruptible
service policies was carried out in \cite{atasol}.
The asymptotics of the queueing control problem were shown to be
governed by a BCP that is a special case of the one identified in this
paper. Thanks to the more special model and the (easier to handle)
interruptible service
assumption, it was possible to attain an analogous result for the
\textit{headcount} process as well
as queue-length (in this paper, Theorem \ref{th1} is proved for the
queuelength process only).

We use the following notation throughout the paper. For a positive
integer $d$ and $x\in\R^d$,
we write $\|x\|$ for $\sum_{l=1}^d|x_l|$, and for $f\dvtx[0,\iy)\to\R^d$,
$\|f\|_t=\sup_{0\leq s\leq t}\|x(s)\|$. We denote by $\mD_{\R^d}$ the
space of functions from $[0,\infty)$ to $\R^d$ that are right
continuous with left limits (RCLL), and equip it with the usual
Skorohod topology. We remove the subscript when $d=1$. For a sequence
of r.v.'s $\{X^n\}, X$, with values in $\R^d$, or processes with sample
paths in $\mD_{\R^d}$, $X^n\Rightarrow X$ denotes convergence in distribution.

The remainder of the paper is organized as follows. Section~\ref{sec:model} describes the model
and states the result regarding the lower bound. It also contains a
discussion of sojourn time costs
as well
as an aspect of the lower bound related to pathwise dominance.
The definitions of various diffusion-scaled processes and some
preliminary lemmas appear in Section~\ref{sec:preliminaries}. This
section also contains a formulation and solution of the underlying
Brownian control problem.
The proof of the lower bound appears in Section~\ref{sec:lowerbound}.
Section~\ref{sec:asymopt} contains the second main result, asserting
that the lower bound is tight. This is shown by constructing a sequence
of policies that asymptotically achieves the bound.

\section{The model and a lower bound on performance} \label{sec:model}

\subsection{Model, scaling, heavy traffic assumptions}\label{sec21}

\begin{figure}[b]

\includegraphics{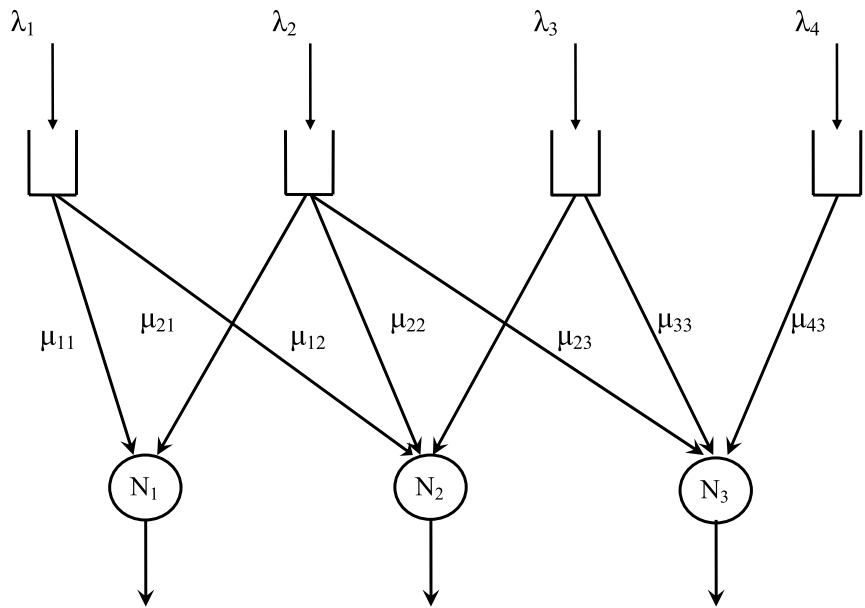}

\caption{A parallel server system.}\label{fig:SBR}
\end{figure}

We consider a sequence, indexed by $n\in\N$, of parallel server
systems with
$\bI$ classes of customers
and $\bJ$ pools of servers, an example of which is depicted in Figure~\ref{fig:SBR}.
The index set for the classes is denoted by $\calI$ and that for
the pools by $\calJ$ (thus $|\calI|=\bI$ and $|\calJ|=\bJ$).

%

Arrivals are modeled as independent renewal processes, denoted by
$(A^n_i,\break  i\in\calI)$ so that $A_i^n(t)$ is the number of class-$i$
arrivals by time $t$. To construct these processes, let $(A_i, i\in
\calI
)$ be independent (undelayed) renewal processes, where,
for each $i$, the time of the first arrival and the inter-arrival times
are positive i.i.d. random
variables with mean 1, variance $(C_{i,{\IA}})^2\ge0$. We assume that
the inter-arrival times have a finite moment of some order $r>2$. The
processes $A^n_i$ are defined as time accelerations of the above, namely
\[
A^n_i(t)=A_i\bigl(\lambda_i^nt
\bigr),\qquad t\ge0, i\in\calI,
\]
where the acceleration parameters satisfy $\lim_n\lambda
_i^n/n=\lambda_i>0$, and
%
%
\begin{equation}
\label{08} \wh\lambda_i^n:=n^{-1/2}\bigl(
\lambda_i^n-n\lambda_i\bigr)\to\wh\lambda
_i\in(-\iy,\iy), \qquad i\in\I
\end{equation}
as $n\to\iy$.

For $j\in\calJ$,
$N_j^n$ denotes the number of servers in pool $j$ and is assumed to satisfy
%
%
\begin{equation}
\label{01} N^n_j=\nu_jn^{1/2}+O
\bigl(n^{1/4}\bigr),\qquad j\in\calJ
\end{equation}
for some constants $\nu_j>0$. We assume that service times are
exponentially distributed and denote by $1/\mu_{ij}^n$ the mean service
time of a class-$i$ customer with a server from pool $j$ (so that $\mu
_{ij}^n$ is interpreted as the corresponding rate). If servers from
pool $j$ cannot serve customers from class $i$, write $\mu_{ij}^n=0$.
This property is assumed to be independent of $n$. Write $i\sim j$
or $j\sim i$ if $\mu_{ij}^n>0$ (for all, equivalently, one $n$). This
information is encoded in
a graph $\calG$ whose vertex set is $\calI\cup\calJ$ and has an edge
connecting $i\in\calI$
and $j\in\calJ$ if and only $i\sim j$. The edge set of the graph is
denoted by $\calE$, and thus $(i,j)\in\calE$
if and only if $i\sim j$. Denote by $\bK$ the cardinality of $\calE$.
When all servers from pool $j$ are occupied with class-$i$ customers,
they jointly serve this class at rate $N^n_j\mu_{ij}^n$. Further assume
that
%
%
\begin{equation}
\mu_{ij}^n=\mu_{ij} n^{1/2}+O
\bigl(n^{1/4}\bigr),
\end{equation}
so that
%
%
\begin{equation}
\label{35} \bar\mu_{ij}^n:=n^{-1}N^n_j
\mu_{ij}^n\to\bar\mu_{ij}:= \mu _{ij}
\nu_j\qquad \mbox{as } n\to\iy,
\end{equation}
and assume that $\bar\mu_{ij}>0$ whenever $i\sim j$ (clearly, $\bar
\mu
_{ij}=0$ otherwise).
Also, assume that, for every $i\in\calI$ and $j\in\calJ$,
%
%
\begin{equation}
\label{36} \wh\mu_{ij}^n:=n^{1/2}\bigl(\bar
\mu_{ij}^n-\bar\mu_{ij}\bigr)\to\wh\mu
_{ij}\in (-\iy,\iy)\qquad \mbox{as } n\to\iy.
\end{equation}
Thus, the nominal joint processing rate of pool-$j$ servers for
class-$i$ customers
(namely $N_j^n\mu_{ij}^n$) is asymptotic to $n\bar\mu_{ij}$. At the
same time, the
rate of an individual
server (namely $\mu_{ij}^n$) is asymptotic to $n^{1/2}\mu_{ij}$. The
quantities $\bar{\mu}=(\bar\mu_{ij})$ will appear in the fluid model,
while $\mu=(\mu_{ij})$ will show in the diffusion model.

Let $\Sig$ denote the set of $\bI\times\bJ$ matrices for which the
$(i,j)$ entry is zero whenever $i\nsim j$.
Let $\X$ be the subset of $\Sig$ of ``column-substochastic'' matrices.
That is,
members $\xi$ of $\X$ satisfy
$\xi_{ij}\ge0$ for every $i,j$, $\sum_i \xi_{ij}\le1$ for every
$j$ and
$\xi_{ij}=0$ for $(i,j)\notin\calE$.
Following \cite{mansto}, for $\xi\in\Sig$, write $\bar\mu(\xi)$
for the
column vector
$(\bar\mu(\xi)_1,\bar\mu(\xi)_2,\ldots,\bar\mu(\xi)_\bI)'$, where
%
%
\begin{equation}
\label{39} \bar\mu(\xi)_i=\sum_{j\in\calJ}\bar
\mu_{ij}\xi_{ij},\qquad i\in\I.
\end{equation}
The first order parameters $\lambda=(\lambda_i)$ and $\bar\mu=(\bar
\mu_{ij})$
are assumed to satisfy a critical load condition. To specify it,
consider a \textit{static fluid model}, consisting of $\bI$ classes of
fluid and $\bJ$ processing stations.
For $i\in\calI$, fluid of class $i$ enters at rate $\lambda_i$.
Each station may divide its processing effort so as to process fluids
of different classes simultaneously.
A member $\xi\in\X$ is said to be an \textit{allocation matrix} for the
model, representing how the effort is distributed among classes.
When the system operates under a given allocation matrix $\xi$, the
element $\xi_{ij}$ represents the fraction of effort devoted at station
$j$ to processing
class-$i$ fluid. Consequently, station $j$ processes class-$i$ fluid at
rate $\bar\mu_{ij}\xi_{ij}$. A system operating under $\xi$ is balanced
if the \textit{balance equation} $\bar\mu(\xi)=\lambda$ holds so that the
stations process the incoming fluid at the rate at which it enters the system.

Consider now the following linear program:
%
%
\begin{equation}
\label{15} \mbox{Minimize $\rho$ over $\xi\in\X$ subject to $\bar\mu(\xi )=
\lambda$ and $\displaystyle\sum_i\xi_{ij}\le\rho$, $j
\in\calJ$, $\rho\geq0$}.\hspace*{-35pt}
\end{equation}
The following condition asserts that the static fluid model is
critically loaded.

\begin{assumption}[{[Heavy traffic (HT)]}]\label{assn2}
There exists a unique optimal solution
$(\xi^*,\rho^*)$ to the linear program \eqref{15}. Moreover,
$\sum_{i\in\calI}\xi_{ij}^*=1$ for all $j\in\calJ$ (and consequently,
$\rho^*=1$).
\end{assumption}
Following terminology from \cite{harlop}, a pair $(i,j)$ where $i\sim
j$, is called an \textit{activity},
and an activity $(i,j)$ is said to be \textit{basic} if $\xi^*_{ij}>0$ and
is said to be \textit{nonbasic} otherwise.
In the static fluid model operating under $\xi^*$, class-$i$ fluid is
processed at a positive rate by station $j$
if and only if $(i,j)$ is a basic activity.
Let $\calG_{b}$ be the sub-graph of $\calG$, having $\calI\cup\calJ
$ as
a vertex set, and an edge connecting
$i\in\calI$ and $j\in\calJ$ if and only if $(i,j)$ is a basic activity.
The edge set of $\calG_b$ will be denoted
by~$\calE_b$. We will write $\calE_{nb}=\calE\setminus\calE_b$ for the
set of nonbasic activities. For $i\in\I$ we let $\J(i):=\{j\in\J\dvtx
(i,j)\in\calE_{b}\}$ be the set of server pools that are connected to
class $i$ via basic activities and, similarly, for $j\in\J$, we let
$\I
(j):=\{i\in\I\dvtx(i,j)\in\calE_{b}\}$ be the set of customer classes
connected to server pool $j$ via basic activities.
%
\begin{assumption}[(Complete resource pooling (CRP) \cite{harlop})]\label{assn3}
The graph $\calG_b$ is connected.
\end{assumption}
Both Assumptions \ref{assn2} and \ref{assn3} will be in force
throughout the paper.
In the context of conventional heavy-traffic, the connectedness of the
stations via basic activities leads to a high level of cooperation
in that the system asymptotically behaves as if it operates under a
single super-server
\cite{belwil2,mansto}. As explained in \cite{harlop} (see also
\cite{belwil1,belwil2,wil}),
the CRP condition is related to the so-called \textit{workload process}
being one-dimensional, and allows for the corresponding Brownian
control problem to admit a one-dimensional solution.
It is known from Williams \cite{wil} that, under Assumption \ref
{assn2}, $\calG_b$ is \textit{connected} if and only if it is a \textit{tree}. Both the lower bound and the asymptotic optimality results build
on the tree structure of $\calG_b$.

It will be useful to state an equivalent form of the above assumptions,
given by Mandelbaum and Stolyar \cite{mansto}.
To this end, denote
\[
\scrM=\bigl\{\bar\mu(\xi)\dvtx\xi\in\X\bigr\}.
\]
Note that $\scrM$ is a convex polygonal domain and a subset of $\R
_+^\bI$.
It is argued in~\cite{mansto} that
the conjunction of the HT and the CRP conditions is equivalent to the
following statement:
\textit{$\lambda$ is a maximal element of $\scrM$ w.r.t. the usual partial
order on $\R_+^\bI$\textup{;}
the unit outward normal to $\scrM$ at $\lambda$ is unique\textup{;} and the matrix
$\xi\in\X$ for which $\bar\mu(\xi)=\lambda$
is unique.} (Note that, because $\lambda_i>0$ for all $i\in\calI$, it
follows that $\scrM$ is an $\bI$-dimensional set,
as assumed a priori in \cite{mansto}.) We will denote the unit outward
normal to $\scrM$ at $\lambda$ by $\theta$. As argued in \cite{mansto},
$\theta$ must satisfy $\theta_i>0$, $i\in\I$. These facts will be used
in our analysis.

We continue the description of the probabilistic model. We let a
complete probability space $(\Om,\calF,\PP)$ be given, supporting
all random variables and stochastic processes defined below.
We write $\EE$ for expectation w.r.t. $\PP$.
Let $B^n_{ij}$ denote the process representing the number of pool-$j$
servers working on class-$i$ customers [note that
$B^n_{ij}\equiv0$ for $(i,j)\notin\calE$].
Let $X^n_i$, $Q^n_i$ and $I^n_j$ denote, respectively, the number of
class-$i$ customers in the system (the ``headcount'' process),
the number of class-$i$ customers in the queue and the number of
pool-$j$ servers that are idle.
Note that
%
%
\begin{equation}
\label{06} X^n_i=Q^n_i+\sum
_jB^n_{ij},\qquad i\in\calI
\end{equation}
and
%
%
\begin{equation}
\label{07} N^n_j=I^n_j+\sum
_iB^n_{ij},\qquad j\in\calJ.
\end{equation}
We are given standard (unit rate) Poisson processes $(S_{ij}$,
$(i,j)\in
\calE)$. The number of service
completions of class-$i$ customers by pool-$j$ servers is constructed
by setting
%
%
\begin{equation}
\label{09} D^n_{ij}(t)=S_{ij}
\bigl(T^n_{ij}(t)\bigr),
\end{equation}
where
%
%
\begin{equation}
\label{10} T^n_{ij}(t)=\mu^n_{ij}\int
_0^tB^n_{ij}(s)\,ds.
\end{equation}
We then have
%
%
\begin{equation}
\label{11} X^n_i(t)=X^n_i(0)+A^n_i(t)-
\sum_jD^n_{ij}(t).
\end{equation}
Naturally, it is required that for all $t\geq0$,
%
%
\begin{equation}
\label{eq:nonneg} X_{i}^n(t),B_{ij}^n(t),Q_i^n(t)
\in\mathbb{Z}_+,\qquad i\in\I,j\in\J.
\end{equation}

For each $n$, the processes $(A_i, i\in\calI)$, $(S_{ij}, (i,j)\in
\calE)$
and the initial condition
$((B^n_{ij}(0),(i,j)\in\calE)), Q^n_i(0))$ are assumed to be mutually
independent.
We refer to $(A,S,B^n(0),Q^n(0))$ as the \textit{primitives}.

When routing decisions are made in a causal manner based on the
observed histories of the processes involved, namely
%
%
\begin{equation}
\label{79} \PI^n:=\bigl(X^n,Q^n,B^n,I^n,D^n,T^n
\bigr),
\end{equation}
the construction of the departure process via \eqref{09} and \eqref
{10} assures
that the customers' service times are independent, exponential random variables
(in particular, this follows from \cite{bremaud}, Theorem 16, page 41;
see also the discussion in Section~2.1 of \cite{DaiTezcanSSC}).
For the treatment of
this paper, it will not be necessary to require any nonanticipating
property of the class of
policies we consider (although the exponential property will be lost
when the policy
does not satisfy a suitable nonanticipating property). Instead, we use an
elaborate definition of the term ``policy'' that merely relies on the
\textit{equations} presented thus far. More precisely, \textit{any process
$\PI^n$, of the form \textup{\eqref{79}} and possessing RCLL sample paths
will be referred to as a policy for the $n$th system, provided that
equations \textup{\eqref{06}}--\textup{\eqref{eq:nonneg}} hold, and
that the stochastic primitives satisfy our probabilistic assumptions
mentioned above.}
Given $n$, the collection of all policies for the $n$th system is
denoted by $\calP^n$. Note that policies need not satisfy any work
conservation condition.

\subsection{Cost functional and asymptotic lower bound}\label{sec22}

Our results will be concerned with asymptotically minimizing a cost
associated to the diffusion-scaled queueing process, defined by
%
%
\begin{equation}
\label{04} \wh Q_i^n(t)=n^{-1/2}Q_i^n(t),\qquad
i\in\calI.
\end{equation}
Let a cost function $C\dvtx\R_+^\bI\to\R_+$ be given, that is continuous
and nondecreasing with respect to the usual partial order on $\R_+^\bI
$. Fix $u>0$. The cost criterion of interest will be
%
%
\begin{equation}
\label{78} \int_0^u C\bigl(
\hQ^n(s)\bigr)\,ds.
\end{equation}
Note that this criterion is a r.v. for each $n$. We do not formulate
the problem in terms of minimizing
the \textit{expected value} of \eqref{78}. Considering \eqref{78} allows
us to state a result on the asymptotic behavior of these r.v.'s that is
stronger than one about their expectation; see Remark \ref{rem2}.
%
\begin{assumption}\label{asum:Cstar}
The function $C_*(\cdot)$ defined by
%
%
\begin{equation}
C_*(a)=\inf\bigl\{C(q)\dvtx q\in\R_+^\bI, \theta'q=a
\bigr\},\qquad a\ge0, \label{eq:cstartdefin}
\end{equation}
is convex.
\end{assumption}
Clearly, a sufficient condition
for the convexity of $C_*$ is the convexity of the function $C$.
However, it is not necessary. Consider, for example, $\I=\{1,2\}$ and
the cost function $C(x_1,x_2)=2(x_1+x_2)^2-(x_1-x_2)^2$ for $x\in\R
_+^2$ and assume $\theta=(1,1)$. Then $C_*(y)=y^2$ is convex while $C$
is not.

To state our first main result we introduce additional notation: for
$x\in\R$, $x^+=\max(x,0)$ and $x^-=\max(-x,0)$. The \textit{Skorohod map}
$\Gam\dvtx\mD([0,\iy)\dvtx\R)\to\mD([0,\iy),\R_+)$ is defined by
%
%
\begin{equation}
\label{54} \Gam[\zeta](t)=\zeta(t)+\sup_{s\le t}\bigl(-\zeta(s)
\bigr)^+, \qquad t\ge0.
\end{equation}

The process
%
%
\begin{equation}
\wh X_i^n(t)=n^{-1/2} \biggl(X_i^n(t)-
\sum_{j}\xi _{ij}^*N_j^n
\biggr) \label{eq:hXdefin}
\end{equation}
represents the diffusion-scale deviation of the headcount process from
the quantities dictated by the static fluid model.
Throughout, we assume
%
%
\begin{equation}
\label{34} \wh X^n(0)\To X_0\qquad\mbox{as } n\to\iy,
\end{equation}
where $X_0$ is a.s. finite r.v.
Finally,
\[
\ell_i:=\wh\lambda_i-\sum
_j\wh\mu_{ij}\xi^*_{ij}\quad\mbox{and}\quad
\s_i:= \biggl(\lambda_iC_{i,\IA}^2+
\sum_j\bar\mu_{ij}\xi^*_{ij}
\biggr)^{1/2},\qquad i\in \calI.
\]

\begin{theorem}\label{th1}
Fix an arbitrary sequence $\{\PI^n=(X^n,Q^n,B^n,I^n,\break D^n,T^n), n\in\N\}$
of policies. Then $\{\PI^n\}$ can be coupled on a common probability
space with the r.v. $X_0$ and an $\bI$-dimensional Brownian motion $W$
(with drift vector $\ell$ and diffusion coefficient $\sigma$) so that
$W$ and $X_0$ are mutually independent and, w.p.1,
\[
\liminf_{n\to\iy}\int_0^u C
\bigl(\wh Q^n(t)\bigr)\,dt \ge\int_0^uC_*
\bigl(Q^*(t)\bigr)\,dt,
\]
where $Q^*$ is the (one-dimensional) reflected Brownian motion given
by\break
$\Gam[\theta'X_0+\theta'W]$.
\end{theorem}

\begin{remark} \label{rem2}
A more standard control theoretic setting is one where an \textit{expected} cost, such as
$\EE[\int_0^u C(\wh Q^n(t))\,dt]$, is to be minimized. An asymptotic
lower bound on the expected cost follows immediately from the above
result, using
Fatou's lemma, namely $\liminf\EE[\int_0^u C(\wh Q^n(t))\,dt] \ge\EE
[\int_0^uC_*(Q^*(t))\,dt]$.
The result stated in the theorem is, of course, stronger.
\end{remark}

\begin{remark}
The family $\mathcal{P}^n$ includes both preemptive and nonpreemptive
policies. The policy that we will construct in Section~\ref{sec:asymopt} is
nonpreemptive, but we will prove that it is asymptotically optimal
within the larger family~$\mathcal{P}^n$.
\end{remark}

\subsection{Discussion}\label{sec23}

\textit{On sojourn time costs.}
In the NDS regime, unlike in the conventional regime,
sojourn time costs lead to control policies that are distinct from
those designed
to minimize waiting times. We provide here a heuristic argument,
showing that sojourn time costs
can be expressed as queue-length (or waiting time) costs. Rigorizing
and extending this viewpoint is subject for future work.
This heuristic argument provides further motivation to study
queue-length costs.

For this discussion we remove the superscript $n$ from the notation.
For $t\ge0$, $i\in\calI$, denote by $\Del_i(t)$ [resp., $\Sig
_i(t)$, $\mathrm{SOJ}_i(t)$]
the properly scaled waiting time (resp., service time, sojourn time)
experienced by the
class-$i$ customer to first arrive after time $t$.
The scaling is obtained by multiplying the original quantities by
$\sqrt n$
(see~\cite{RamiNDS}).
We have
\[
\mathrm{SOJ}_i(t)=\Del_i(t)+\Sig_i(t).
\]
Consider a cost of the form
\[
J=\int_0^u\sum_iC_i
\bigl(\mathrm{SOJ}_i(t)\bigr)\,dt.
\]

One expects that under reasonable policies the fraction of class $i$
customers routed to servers in pool $j$ be roughly $p_{ij}=\bar{\mu
}_{ij}\xi_{ij}^*/\lambda_i$.
Moreover, one expects that $\Sig_i(t)$ is approximately a mixture of
exponentials where it is an exponential with rate $\mu_{ij}$ with
probability $p_{ij}$ (this is consistent with
the result for the so-called inverted V model in \cite{RamiNDS}; see
Theorem 2.2 there).
Denote by $F_{\Sig_i}$ the corresponding distribution function.
We have $\Ex[C_i(\mathrm{SOJ}_i(t))]=E[R_i(t)]$, where
\[
R_i(t)=\Ex\bigl[C_i\bigl(\Del_i(t)+
\Sig_i(t)\bigr)|\Del_i(t)\bigr].
\]
By the above discussion one expects that $R_i(t)\approx\tilde C_i(\Del_i(t))$,
where
\[
\tilde C_i(y)=\int_0^\iy
C_i(y+x)\,dF_{\Sig_i}(x),\qquad  y\ge0.
\]
Using $\Del_i(t)\approx\lambda_iQ_i(t)$ (by the \textit{snapshot principle})
thus shows
\[
J\approx\Ex \biggl[\int_0^u\sum
_i \tilde C_i\bigl(\lambda_iQ_i(t)
\bigr)\,dt \biggr].
\]

\textit{On pathwise lower bounds}.
Note that Theorem \ref{th1} does not assert that, under the coupling,
w.p.1,
%
%
\begin{equation}
\label{100} \liminf_{n\tinf}C\bigl(\wh{Q}^n(t)\bigr)
\geq C_*\bigl(Q^*(t)\bigr),\qquad t\ge0.
\end{equation}
It only provides an integral version of this inequality.
In conventional HT, \eqref{100} is often called a pathwise lower bound,
and specifically, for the
parallel server model, it is shown to hold in \cite{mansto}.

However, in the NDS regime, particularly, under the setting of Theorem
\ref{th1}, \eqref{100} is a false\vadjust{\goodbreak}
statement. In fact,
\textit{given suitable initial conditions \textup{(}e.g.,~zero\textup{)} one can find
constants $t_0>0$ and $c>0$ such that,
under a suitable policy, $\hat Q^n(t)=0$ for all $t\in
[t_0,t_0+cn^{-1/2}]$, with probability tending to one.}
[This clearly shows that under no coupling can \eqref{100} hold.]
We do not prove this statement here. A detailed study of an analogous phenomenon
in the Halfin--Whitt regime, referred to as \textit{null-controllability},
has been studied in detail
\cite{AtarNC,NC2}.
Briefly, this phenomenon is described as follows. Under suitable
algebraic conditions on the system parameters,
the \textit{critically loaded} parallel server model in the Halfin--Whitt
regime can be controlled so that
all queue-lengths vanish for $O(1)$ units of time, with probability
tending to one.
When the mechanisms described in these works are applied
in the NDS regime, they yield vanishing of all queue-lengths, though
only for $O(n^{-1/2})$ units of time.
[Clearly, such a property could not hold for $O(1)$ units of time in
the NDS regime, since this would
contradict Theorem \ref{th1}.]

This complication explains why the proof of Theorem \ref{th1} is more
involved than
the analogous lower bound in
the conventional regime (as well as the need for tools such as
Proposition \ref{prop2}
and Lemma \ref{lem5}).

\textit{On other heavy traffic regimes.}
Atar \cite{RamiNDS} emphasizes the viewpoint that there exists
a whole spectrum of heavy traffic diffusion approximations.
For any $\al\in[0,1]$ one obtains a distinct heavy traffic regime by setting
the quantities $\lambda$, $\mu$ and $N$ proportional to $n$,
$n^{1-\al}$
and, respectively,
$n^\al$, while maintaining the critical load condition of having
$\lambda-N\mu$ at the order of $n^{1/2}$.
The conventional, NDS and Halfin--Whitt regimes can then be identified
with the cases
$\al=0,1/2$ and, respectively, $1$. Whether results analogous to those
of this paper
hold for all values of $\alpha\in(0,1)$ is left as an important open
question. In view of the so-called null-controllability results in \cite{AtarNC,NC2}, $\al=1$
should be
excluded since the lower bound is not expected to hold in the
Halfin--Whitt regime.

\section{Preliminaries}\label{sec:preliminaries}

\subsection{Diffusion-scale processes}\label{sec:diff_scale}

In the present subsection we define some diffusion-scale processes and
develop relations that they satisfy. Let
%
%
\begin{eqnarray}
\label{45} \bar T^n_{ij}(t)&=&n^{-1}T^n_{ij}(t)
\equiv\frac{\mu^n_{ij}}{n}\int_0^tB_{ij}^n(s)\,ds,
\\
\label{20} \wh A^n_i(t)&=&n^{-1/2}
\bigl(A_i^n(t)-\lambda_i^nt
\bigr),\qquad i\in\calI,
\\
\label{33} \wh S^n_{ij}(t)&=&n^{-1/2}
\bigl(S_{ij}(nt)-nt\bigr) \qquad (i,j)\in\calE,
\\
\label{03} \til B^n_{ij}&=&B^n_{ij}-
\xi^*_{ij}N^n_j,\qquad \wh B^n_{ij}=n^{-1/2}
\til B^n_{ij},
\\
\label{21} V^n_{ij}&=&n^{-1/2}\bigl(D^n_{ij}-T^n_{ij}
\bigr)\equiv\wh S^n_{ij}\circ\bar T^n_{ij},
\\
\label{24} \ell_i^n&=&\wh\lambda_i^n-
\sum_j\wh\mu^n_{ij}
\xi_{ij}^*
\end{eqnarray}
and
%
%
\begin{equation}
\label{18} W^n_i(t)=\ell^n_it+
\wh A_i^n(t)-\sum_j
V_{ij}^n(t).
\end{equation}
Since $\sum_i\xi^*_{ij}=1$, we have by \eqref{07} that
%
%
\begin{equation}
\label{12} I^n_j+\sum_i
\til B^n_{ij}=0,\qquad j\in\calJ,
\end{equation}
and by \eqref{06} and \eqref{eq:hXdefin} that
%
%
\begin{equation}
\label{13} \wh X^n_i=\wh Q^n_i+
\sum_j\wh B^n_{ij}.
\end{equation}
Using \eqref{09}, \eqref{10}, \eqref{11} and \eqref{18}, we get
\begin{eqnarray*}
\wh X_i^n(t) &=& \wh X_i^n(0)+n^{-1/2}A^n_i-n^{-1/2}
\sum_jD_{ij}^n(t)
\\
&=& \wh X_i^n(0)+W^n_i(t)+n^{1/2}
\lambda_i t-n^{1/2}\sum_j
\bar\mu _{ij}\xi ^*_{ij}t\\
&&{}-n^{-1/2}\sum
_j\mu^n_{ij}\int_0^t
\til B^n_{ij}(s)\,ds.
\end{eqnarray*}
Since $\bar\mu(\xi^*)=\lambda$, the third and forth terms above
cancel out
so that letting
%
%
\begin{equation}
\label{eq:eps_defin} \eps^n_{ij}:=n^{-1/2}
\mu^n_{ij}-\mu_{ij},
\end{equation}
we arrive at the following identity:
%
%
\begin{equation}
\label{19} \wh X^n_i(t)=\wh X^n_i(0)+W^n_i(t)-
\sum_j\bigl(\mu_{ij}+
\eps_{ij}^n\bigr)\int_0^t
\til B^n_{ij}(s)\,ds.
\end{equation}
Using \eqref{13}, we obtain
%
%
\begin{equation}
\label{31} \wh Q^n_i(t)=\wh X^n_i(0)+W^n_i(t)-
\sum_j\bigl(\mu_{ij}+
\eps_{ij}^n\bigr)\int_0^t
\til B^n_{ij}(s)\,ds-\sum_j\wh
B^n_{ij}(t).
\end{equation}
Identities \eqref{19} and \eqref{31} will be used in the sequel.

\subsection{Auxiliary results}

\begin{lemma}
\label{lem2}
\textup{(i)} The rescaled primitive processes $(\wh A^n_i,i\in\calI)$ and
$(\wh S^n_{ij},\break (i,j)\in\calE)$ and initial condition $\wh X^n(0)$,
jointly converge in law, uniformly on compacts,
to processes denoted $(W_{A,i},i\in\calI)$ and $(W_{S_{ij}},(i,j)\in
\calE)$, and the r.v. $X_0$,
where $W_{A,i}$ (resp., $W_{S_{ij}}$) is
a zero mean Brownian motion with diffusion coefficient $\lambda
_i^{1/2}C_{i,\IA}$ (resp., 1). Moreover, the $\bI+\bK$ Brownian
motions and
the r.v. $X_0$ are mutually independent.
\begin{longlist}[(iii)]
\item[(ii)] The parameters defined in \eqref{24} and \eqref
{eq:eps_defin} satisfy
%
%
\begin{equation}
\label{44} \ell^n_i\to\ell_i:=\wh
\lambda_i-\sum_j\wh\mu_{ij}
\xi ^*_{ij}\quad \mbox{and}\quad \eps^n_{ij}=O
\bigl(n^{-1/4}\bigr)\qquad \mbox{as } n\to\iy.
\end{equation}
\item[(iii)]
Consequently, the processes
%
%
\begin{equation}
\label{46} \wh W^n_i(t):=\ell^n_it+
\wh A^n_i(t)-\sum_j\wh
S^n_{ij}\bigl(\bar\mu _{ij}\xi ^*_{ij}t
\bigr),\qquad i\in\calI,
\end{equation}
and the initial condition $X^n_0$
jointly converge in law to mutually independent processes $(W_i, i\in
\calI)$ and r.v. $X_0$, where $W_i$ is a Brownian motion
starting from zero, with drift $\ell_i$ [cf. \eqref{44}] and diffusion
coefficient
%
%
\begin{equation}
\label{48} \s_i:= \biggl(\lambda_iC_{i,\IA}^2+
\sum_j\bar\mu_{ij}\xi ^*_{ij}
\biggr)^{1/2},\qquad i\in\calI.
\end{equation}
\end{longlist}
\end{lemma}
Invoking the Skorohod representation theorem, we can, and will, assume
throughout that the
convergence statements of the above lemma occur in an a.s. sense.
%
\begin{remark}
\label{rem1}
Note that Lemma \ref{lem2} deals with convergence of processes that
depend only on the primitives, and
thus the same a.s. limit is attained under any policy.
\end{remark}

\begin{pf*}{Proof of Lemma \ref{lem2}}
(i) It is well known that a renewal process, scaled in the fashion of~\eqref{20}
and \eqref{33},
converges in law, uniformly on compacts, to a Brownian motion with zero
mean and diffusion coefficient as stated
\cite{bil}, Section~17. The mutual independence of the processes and
the independence from the initial conditions
follows the validity of this property for the pre-limit objects. (ii)
The first statement follows by \eqref{08}, \eqref{36} and \eqref{24}.
The second follows by \eqref{01}, \eqref{36} and \eqref{eq:eps_defin}.
\end{pf*}

\subsection{Diffusion model formulations}

In this subsection we present two diffusion models, originating from
\cite{harlop} and \cite{mansto}. We associate to these models a
control problem
analogous to the one used above for the queuing network, and provide a
complete solution.
The optimum is exactly analogous to the lower bound from Theorem \ref{th1}.
The main point of this analysis is that the problem of identifying a
minimizing control
(Proposition \ref{prop1} below)
can later be mimicked to construct a policy for the queueing network
that achieves
the lower bound in an asymptotic sense; see Section~\ref{sec:asymopt}.
Along the way we establish equivalence of the two diffusion models.

Harrison and L\'opez \cite{harlop} present and analyze
a model of controlled diffusion, which stands for
the formal limit of diffusion-scaled processes associated with the
queueing model (in the conventional HT regime).
The diffusion model was later used by Bell and Williams \cite{belwil2}
as the basis for the construction of asymptotically optimal policies
for the queueing model.
This diffusion model, to which we refer as \textit{Model} I, consists of
the r.v.'s $X_{0,i}$ and BMs $W_i$ alluded to in Section~\ref{sec:model},
and in addition, processes $(X_i, i\in\calI)$, $(I_j, j\in\calJ)$ and
$(Y_{ij}, (i,j)\in\calE)$, possessing RCLL sample paths,
which satisfy, in addition, the following relations:
%
%
\begin{eqnarray}
\label{22}&\displaystyle X_i(t)=X_{0,i}+W_i(t)+\sum
_{j:i\sim j}\mu_{ij}Y_{ij}(t)\ge0,\qquad t
\ge0, i\in\calI,&
\\
\label{23} &\displaystyle I_j:=\sum_iY_{ij}
\mbox{ is nondecreasing and $I_j(0)\ge0$,} \qquad j\in\calJ,&
\\
\label{27}&\displaystyle \mbox{$Y_{ij}$ is nonincreasing and $Y_{ij}
\le0$}\qquad (i,j)\in \calE_{nb}.&
\end{eqnarray}
($Y_{ij}$ are further required in \cite{harlop} to be adapted, but here
we take the viewpoint of
\cite{belwil2} where this requirement is dropped.)
Taking formal limits of the scaled processes from our queueing model
gives rise to the same diffusion model.
Indeed, assuming that $W^n$ converges to $W$ and that $\wh B^n$
converge to zero, noting that
$\wh X^n(0)$ converge to $X_0$ and $\eps^n\to0$, and writing $Y$ for a
limit of $-\int_0^\cdot\til B^n$,
equation \eqref{22} arises as a limit form of \eqref{19} [also of
\eqref
{31}], where $X$
stands for the limit of $\wh X^n$ (also of $\wh Q^n$). The
nonnegativity constraint on $X_i$
is clear from that of $Q^n$, observing \eqref{04}. Similarly, if $I_j$
corresponds
to a limit of $\int_0^\cdot I^n_j$, then $I_j$ should be nondecreasing
and, because of
\eqref{12}, satisfy \eqref{23}. Finally, \eqref{27} represents the fact
that $\til B^n_{ij}\ge0$
for any nonbasic activity $(i,j)$ due to \eqref{03} and the fact that
$\xi^*_{ij}=0$ for such an activity.
Making the above formal statements rigorous will be one of the main
issues dealt with in
Sections~\ref{sec:lowerbound} and in~\ref{sec:asymopt}.

Mandelbaum and Stolyar \cite{mansto} construct asymptotically optimal
control policies for
the parallel server model in the conventional HT regime, without
explicitly alluding to a diffusion model (or a Brownian control problem).
However, their verbal discussion and mathematical treatment
of the diffusion scaled processes suggest the following diffusion
model, to which we refer here as \textit{Model} II.
In addition to the random vector $X_0$ and the $\bI$-dimensional
process $W$ from Section~\ref{sec:model},
the model consists of $\bI$-dimensional processes $X$ and $Z$. These
are assumed to have
RCLL sample paths
and satisfy
%
%
\begin{eqnarray}
\label{29}&\displaystyle X_i(t)=X_{0,i}+W_i(t)+Z_i(t)
\geq0, \qquad t\ge0, i\in\I,&
\\
\label{30} &\displaystyle \theta'Z \mbox{ is nondecreasing and $
\theta'Z(0)\ge0$.}&
\end{eqnarray}

The interpretation of $X$ is the same as in Model I, while $Z$
corresponds to the term
$\sum_j\mu_{ij}Y_{ij}$ of \eqref{22}. Below we claim that Models I and~II are equivalent in a suitable sense,
and that they both achieve the lower bound of Theorem \ref{th1}.
%
\begin{proposition} \label{prop1} Suppose that Assumptions \ref{assn2}
and \ref{assn3} hold and fix an initial state $X_0$ and a Brownian
motion $W$ as in Lemma \ref{lem2}.
\begin{longlist}[(iii)]
\item[(i)] Given a pair $(X,Y)$ that satisfies \eqref{22}--\eqref{27},
there exists $Z$ such that the pair $(X,Z)$ satisfies \eqref
{29}--\eqref{30}.
\item[(ii)] Given a pair $(X,Z)$ that satisfies \eqref{29}--\eqref{30},
there exists $Y$ such that $(X,Y)$ satisfies \eqref{22}--\eqref{27}.
\item[(iii)] Let $(X,Y,Z)$ be such that $(X,Y)$ satisfy \eqref
{22}--\eqref{27} and $(X,Z)$ satisfy \eqref{29}--\eqref{30}. Then, with
probability 1, $C(X(t))\geq C_*(Q^*(t))$ for all $t$, where $Q^*$ is as
in Theorem \ref{th1}. Moreover, the lower bound is attainable: there
exist stochastic processes $(X,Y,Z)$ such that $(X,Y)$ [resp., $(X,Z)$]
satisfies \eqref{22}--\eqref{27} [resp. \eqref{29}--\eqref{30}] and
with probability 1,
\[
C\bigl(X(t)\bigr)=C_*\bigl(Q^*(t)\bigr),\qquad t\ge0.
\]
\end{longlist}
\end{proposition}
Below, the notation $\Sig$, $\X$ is as in Section~\ref{sec:model} and,
for $\xi\in\Sig$, $\bar\mu(\xi)$ is as defined in~\eqref{39}.

\begin{lemma}
\label{lem1}
Let $x\in\Sig$ be such that $x_{ij}\le0$ for all $(i,j)\in\calE
_{nb}$, and
$\sum_ix_{ij}\ge0$ for all $j\in\calJ$. Then $\theta'\bar\mu
(x)\ge0$.
Also, if $x_{ij}=0$ for all $(i,j)\in\calE_{nb}$ and $\sum_ix_{ij}=0$
for all $j\in\calJ$ then $\theta'\bar\mu(x)=0$.
\end{lemma}
\begin{pf}
Let us first show that $\xi:=\xi^*-\eps x$ is an element of $\X$,
provided that
$\eps>0$ is sufficiently small. For $(i,j)\in\calE_{nb}$, $\xi
^*_{ij}=0$, and so $\xi_{ij}\ge0$
by the assumptions of the lemma. For $(i,j)\in\calE_b$, $\xi
^*_{ij}>0$, hence
$\xi_{ij}\ge0$ for all sufficiently small $\eps>0$.
Finally, since we assumed that $\sum_ix_{ij}\ge0$ for all $j\in\calJ$
and since $\xi^*\in\X$, we have that
\[
\sum_i\xi_{ij}\le\sum
_i\xi^*_{ij}\le1, \qquad j\in\calJ,
\]
so that $\xi\in\X$. Next, since $\theta$ is an outward normal to
the convex set $\scrM$ at $\lambda$, we have that $\theta'(m-\lambda
)\leq0$
for every $m\in\scrM$.
Since $\bar\mu(\xi^*)=\lambda$ and $\xi\in\X$, we have
%
%
\begin{equation}
\theta'\bar\mu(\eps x) = \theta'\bigl(\bar\mu\bigl(
\xi^*\bigr)-\bar\mu(\xi )\bigr)=\theta '\bigl(\lambda-\bar\mu(\xi)
\bigr)\geq0. \label{eq:vectrineq}
\end{equation}
Since $\eps>0$, the first claim follows.

For the second part, following \cite{mansto}, we claim that there exist
constants $z_j^*, j\in\J$, such that $\theta_i\bar{\mu}_{ij}=z_j^*$ for
all $(i,j)\in\calE_b$. Indeed, by \eqref{eq:vectrineq} $\theta
'\lambda
=\theta'\bar{\mu}(\xi^*)=\sup_{\xi\in\X} \theta'\bar{\mu
}(\xi)$ and, in
turn, we must have that $\theta_i\bar{\mu}_{ij}=\max_{k}\theta
_k\bar{\mu
}_{kj}$ for all $(i,j)\in\calE_b$. Define
%
%
\begin{equation}
z_j^*=\max_{k}\theta _k\bar{
\mu}_{kj}\label{eq:zdefin}.
\end{equation}
Thus, $\theta_i\mu_{ij}=z_j^*$
for all $(i,j)\in\calE_b$. If $x$ is such that $x_{ij}=0$ for all
$(i,j)\in\calE_{nb}$ and $\sum_{i}\xi_{ij}=1$ for all $j\in\J$,
then\vadjust{\goodbreak}
$\xi:=\xi^*-\varepsilon x$ satisfies $\xi_{ij}=0$ for all $(i,j)\in
\calE
_{nb}$ and $\sum_{i}\xi_{ij}=1$ for all $j\in\J$ and, in turn,
\begin{eqnarray*}
\theta'\bar{\mu}(\varepsilon x)&=&\theta'\bigl(\bar{
\mu}\bigl(\xi^*\bigr)-\bar{\mu }(\xi )\bigr)=\sum_{(i,j)\in\calE_{b}}
\theta_i\mu_{ij}\bigl(\xi_{ij}^*-\xi
_{ij}\bigr)\\
&=&\sum_jz_j^*
\biggl(\sum_{i}\xi_{ij}^*-\sum
_{i}\xi_{ij}\biggr)=0.
\end{eqnarray*}
\upqed\end{pf}

\begin{pf*}{Proof of Proposition \ref{prop1}}
(i) Suppose that relations \eqref{22}--\eqref{27} hold. Letting
$Z_i=\sum_j\mu_{ij}Y_{ij}$, the relation
\eqref{29} is immediate. To show \eqref{30}, we must prove that, for
$0\le s\le t$,
\[
\sum_{ij}\theta_i\mu_{ij}
\bigl(Y_{ij}(t)-Y_{ij}(s)\bigr)\ge0\qquad\biggl[\mbox{resp.},
\sum_{ij}\theta_i\mu_{ij}Y_{ij}(t)
\ge0  \biggr],
\]
which can alternatively be stated as $\theta'\bar\mu(x)\ge0$, where
$x_{ij}=\nu_j^{-1}(Y_{ij}(t)-Y_{ij}(s))$ [resp., $x_{ij}=\nu
_j^{-1}Y_{ij}(t)$].
Properties \eqref{23} and \eqref{27} of $Y$ and an application of Lemma~\ref{lem1} imply $\theta'\bar\mu(x)\ge0$,
in both cases, and the claim is proved.

(ii) Assume we are given a solution $(X,Z)$ to \eqref{29} and \eqref
{30}. We will construct a process $Y$ so that the pair $(X,Y)$
satisfies \eqref{22}--\eqref{27}. Fix an arbitrary $l\in\J$ and set
$I_j= 0$ for all $j\neq l$ and $I_{l}(t)=\frac{\nu_l}{z_l^*}\theta
'Z(t)$ [where $z^*$ is as in \eqref{eq:zdefin}]. Then, given $X_0$,
$X$, $Z$ and $W$, let us show there exists a solution $Y$ to the set of
equations
%
%
\begin{eqnarray}
\sum_{j:(i,j)\in\calE_b}\mu _{ij}Y_{ij}&=&X_i-X_{i,0}-W_i,\qquad
i\in\I,\label{eq:Yconstruct1}
\\
\sum_{i:(i,j)\in\calE_b}Y_{ij}&=&0\qquad \forall j\in\calJ
\setminus\{ l\}, \label{eq:Yconstruct2}
\\
\sum_{i:(i,l)\in\calE_b}Y_{il}&=&\frac{\nu_l}{z_l^*}\theta
'Z,\label {eq:Yconstruct3}
\end{eqnarray}
and $Y_{ij}= 0$ for $(i,j)\notin\calE_b$.
To this end, we first remove one on the $I+J$ equations from the
system. Fix $i_0\in\calI$ and assume that \eqref{eq:Yconstruct1} holds
for all $i\neq i_0$. Then
%
%
\begin{eqnarray}\label{eq:reduction}
\sum_{i\neq i_0}\theta_iX_i(t)&=&
\sum_{i\neq
i_0}\theta _iX_{i,0}+\sum
_{i\neq i_0}\theta_iW_i(t)+\sum
_{i\neq i_0}\theta _i\sum
_{j}\bar{\mu}_{ij}Y_{ij}/
\nu_j
 \nonumber
 \\[-8pt]
 \\[-8pt]
 \nonumber
&=&\sum_{i\neq i_0}\theta_iX_{i,0}+
\sum_{i\neq i_0}\theta _iW_i(t)-
\sum_{j}\frac{\bar{z}_j^*}{\nu_j}Y_{i_0j}+
\theta'Z(t),
\end{eqnarray}
where we used $\theta_i\bar{\mu}_{ij}=z_j^*$ for all $(i,j)\in\calE_b$
(see the proof of Lemma \ref{lem1}) and equations \eqref
{eq:Yconstruct2} and \eqref{eq:Yconstruct3}. On the other hand, since
$(X,Z)$ solves \eqref{29} and \eqref{30} we have that $\theta
'X(t)=\theta'X_0+\theta'W(t)+\theta'Z(t)$ and, together with \eqref
{eq:reduction}, that
\[
\theta_{i_0}X_{i_0}(t)=\theta_{i_0}X_{i_0}+
\theta_{i_0}W_{i_0}(t) +\sum_{j}
\frac{\bar{z}_j^*}{\nu_j}Y_{i_0j}.
\]
Substituting $\theta_i\bar{\mu}_{ij}=z_j^*$ we have thus shown that
\eqref{eq:Yconstruct1} holds for $i_0$ provided that it holds for all
$i\neq i_0$ and provided that \eqref{eq:Yconstruct2} and \eqref
{eq:Yconstruct3} hold. Hence we can remove the equation for $i=i_0$
from \eqref{eq:Yconstruct1}.

The reduced system of equations has $I+J-1$ variables (one for each
basic activity) and the same
number of equations.
By the discussion on page 348 of \cite{harlop}, these equations are
linearly independent.
As a result, given $X_0$, $X$, $Z$ and $W$, there exists a unique
solution $Y$.

Note that \eqref{23} and \eqref{27} hold by construction.

To establish the lower bound in item (iii), let $M(t):=\sup_{s\leq
t}(\theta' X_{0}+\theta' W(t))^{-}$ and note that, by the minimality of
Skorohod problem (see, e.g., \cite{cheman}, Section~2)
\[
\theta' X(t)\geq\theta' X_0+
\theta' W(t)+M(t)\quad\mbox{and}\quad \theta' \tilde{X}(t)\geq
\theta' X_0+ \theta' W(t)+M(t).
\]
Since, $C(x)\geq C_*(\theta'x)$ for all $x\in\mathbb{R}^I$, the lower
bound is established.

Finally, to show that the lower bound is attained we explicitly
construct a pair $(X,Z)$ that attains the lower bound. A process
$(X,Y)$ that attains the lower bound will then be constructed from
$(X,Z)$ as above. To that end, for $t\geq0$, let $Q^*(t)=\theta
'X_0+\theta'W(t)+M(t)$ where $M(t)$ is as above. Set
$X_i(t)=f_i^*(Q^*(t))$ where, given $a\in\mathbb{R}_+$,
$f^*(a):=(f_1^*(a),\ldots,f_I^*(a))$ satisfies
\[
f^*(a)\in\mathop{\argmin}_{q\geq0}\bigl\{C(q)\dvtx q\in\R_+^\bI,
\theta'q=a\bigr\}.
\]
$f^*$ can be selected to be measurable, as follows from Corollary 10.3
in the Appendix of \cite{ethier1986markov}, using the continuity of
$C(\cdot)$. Thus
$\theta'X(t)=Q^*(t)$ and $C(X(t))=C(f^*(\theta
'X(t)))=C(f^*(Q^*(t)))=C_*(Q^*(t))$.
Setting $Z_i:=X_i-X_{0,i}-W_i$, we have that $\theta'Z(t)=\theta
'X-\theta'X_{0}-\theta'W=Q^*(t)-\theta'X_0-\theta'W=M(t)$ so that
$\theta'Z$ is nonnegative and nondecreasing and the pair $(X,Z)$
satisfies \eqref{29} and attains the lower bound.
\end{pf*}

\section{Proof of the lower bound}\label{sec:lowerbound}

In this section we prove Theorem \ref{th1}. The main estimate on which
the proof is based, Proposition \ref{prop2}, is stated in Section~\ref{sec41} where it is also used to prove Theorem \ref{th1}. Proposition
\ref{prop2} is then proved in Section~\ref{sec42}.

\subsection{\texorpdfstring{Proof of Theorem \protect\ref{th1}}
{Proof of Theorem 2.1}}\label{sec41}

An outline of the proof is as follows.
Fix $u>0$, which will serve as a time horizon. Given a sequence of
policies $(\PI^n, n\in\N)$ we show that up to a certain random time
$\tau_n\wedge u$, the cumulative process $\int_0^t\theta' \wh Q^n(s)\,ds$
is asymptotically bounded from below by the integrated RBM $\int_0^t
Q^*(s)\,ds$. We then show that if $\tau_n<u$, then $\theta' \wh Q^n$ is
large on a subinterval of $(\tau_n,u]$ and thus, with high probability,
is bounded from below by the integrated RBM. The convexity of
$C_*(\cdot
)$ is then used in translating these bounds to bounds on the cost.

We turn to the proof. Denote $\eps_M^n=\max_{i,j}|\eps^n_{ij}|$, and
recall that $\eps^n_M=O(n^{-1/4})$ (Lemma \ref{lem2}).
Let us fix a sequence $\varrho_n$ such that
%
%
\begin{equation}
\label{51-} n^{-1/8}\varrho_n\to\iy\qquad \mbox{while }
n^{-1/4}\varrho_n\to0.
\end{equation}
In particular, $\varrho_n$ satisfies
%
%
\begin{equation}
\label{51} n^{-1/2}\varrho_n\to0,\qquad \eps_M^n
\varrho_n\to0,\qquad n^{-1/2}\bigl(\eps^n_M
\bigr)^{-1}\varrho_n^2\to\iy.
\end{equation}
Let a sequence of policies $(\PI^n, n\in\N)$ be given, and define
%
%
\begin{equation}
\label{32} \tau_n=\inf \biggl\{t\ge0\dvtx\max_{i,j}
\biggl|\int_0^t\til B_{ij}^n(s)\,ds
\biggr|\ge\varrho_n \biggr\}\w u.
\end{equation}
Since the conclusion of Theorem \ref{th1} clearly holds when the
function $C_*$ is constant,
we henceforth assume that $C_*$ is not constant. Denote $E^n(t)=\theta
'\wh Q^n(t)$, $t\ge0$. Below $Q^*$ is as defined in Theorem \ref{th1}.
%
\begin{proposition}
\label{prop2}
There exist constants $c$, $\til c$, $\bar c$, a strictly positive
sequence $\{t_n\}$ satisfying
$t_n\varrho_n\to\iy$, a sequence of events $\Om^n$ satisfying
$1_{\Om
^n}\to1$ a.s.,
and processes $P^{*,n}$ and $H^{*,n}$, such that, with
\[
\til\tau_n:=(\tau_n+t_n)\w u,
\]
the following statements hold:
\begin{longlist}[(iii)]
\item[(i)]
$P^{*,n}$ converges to $Q^*$ uniformly on $[0,u]$, a.s.

\item[(ii)] $|H^{*,n}(t)|\le\bar c$, for every $n$ and $t\in[0,u]$, and
%
%
\begin{equation}
\label{70-} \biggl|\int_0^tH^{*,n}(s)\,ds \biggr|
\le\til cn^{-1/2}\varrho_n,\qquad t\in[0,u].
\end{equation}
\item[(iii)] On $\Om^n$ one has
%
%
\begin{eqnarray}
E^n(t) &\ge& P^{*,n}(t)+ H^{*,n}(t)\qquad  \mbox{for
all } t\in[0,\tau_n); \label{70}
\\
E^n(t) &\ge& c\varrho_n \qquad \mbox{for all } t\in[
\tau_n,\til\tau_n), \mbox{ whenever }
\tau_n<u; \label{71}
\\
E^n(t) &\ge&0 \qquad \mbox{for all } t\in[\til\tau_n,u],
\mbox{ whenever } \til\tau_n<u. \label{72}
\end{eqnarray}
\end{longlist}
\end{proposition}
The proof is deferred to Section~\ref{sec42}.
Theorem \ref{th1} will be deduced from the Proposition \ref{prop2},
with the aid of the following lemma for which we define for $-\infty<
a\leq b <\infty$,
for $x\dvtx[a,b]\to\R$ and $\delta>0$,
%
%
\begin{equation}
\bar w_{[a,b]}(x,\delta)=\sup_{s,t\in[a,b];|s-t|\le\delta}\bigl|x(s)-x(t)\bigr|.
\label{eq:wdefin}
\end{equation}
%
\begin{lemma}
\label{lem5}
Let $C_1\dvtx\R\to\R_+$ be a nondecreasing, convex function. Let $T>0$,
$0<\Delta<T/2$, $r>0$, and functions
$p,h\dvtx\R_+\to\R$ be given such that
\[
\biggl|\int_0^th(s)\,ds\biggr |\le\eps, \qquad t\in[0,T]
\]
and
\[
\bigl|p(t)\bigr|+\bigl|h(t)\bigr|\le r,\qquad t\in[0,T].
\]
Then
\[
\int_0^TC_1\bigl(p(t)+h(t)
\bigr)\,dt\ge\int_0^TC_1\bigl(p(t)
\bigr)\,dt-\gamma_1T-\gamma _2T-\gamma_3,
\]
where
\[
\gamma_1=\bar w_{[-r,r]} \biggl(C_1,
\frac{2\eps}{\Del} \biggr),\qquad \gamma _2=\bar w_{[-r,r]}
\bigl(C_1, \bar w_T(p,\Del) \bigr),\qquad \gamma_3=2C_1(r)
\Del.
\]
\end{lemma}
\begin{pf}
For $t\in[0,T-\Del]$, using Jensen's inequality,
\begin{eqnarray*}
\frac1\Del\int_t^{t+\Del}C_1
\bigl(p(t)+h(t)\bigr)\,dt &\ge &C_1 \biggl(\frac1\Del\int
_t^{t+\Del}\bigl(p(s)+h(s)\bigr)\,ds \biggr)
\\
&\ge &C_1 \biggl(\frac1\Del\int_t^{t+\Del}p(s)\,ds
\biggr)-\gamma_1
\\
&\ge& C_1\bigl(p(t)\bigr)-\gamma_1-
\gamma_2.
\end{eqnarray*}
Thus
\[
\int_0^TC_1\bigl(p(s)+h(s)
\bigr)\,ds\ge\int_\Del^{T-\Del}C_1\bigl(p(t)
\bigr)\,dt-\gamma _1T-\gamma_2T.
\]
The result follows.
\end{pf}

\begin{pf*}{Proof of Theorem \ref{th1}}
Extend $C^*$ to $\R$ by setting $C^*=C^*(0)$ on $(-\iy,0)$.
Fix $\Del\in(0,u/2)$.
Let $P^{*,n}$, $H^{*,n}$, $\Om^n$, $\til\tau_n$ be as in Proposition
\ref{prop2}.
Combining Proposition \ref{prop2} and Lemma \ref{lem5} we have, on
$\Om^n$,
\begin{eqnarray*}
\int_0^uC_*\bigl(E^n(t)\bigr)\,dt &
\ge&\int_0^{\tau_n}C_*\bigl(P^{*,n}(t)+H^{*,n}(t)
\bigr)\,dt+C_*(c\varrho _n) (\til\tau _n-\tau_n)
\\
&\ge&\int_0^{\tau_n}C_*\bigl(P^{*,n}(t)
\bigr)\,dt+C_*(c\varrho_n) (\til\tau _n-\tau
_n)-\gamma^n_1u-\gamma^n_2u-
\gamma^n_3,
\end{eqnarray*}
where
\begin{eqnarray*}
\gamma^n_1&=&\bar w_{[-r_n,r_n]} \biggl(C_*,
\frac{2\til cn^{-1/2}\varrho
_n}{\Del} \biggr),\qquad \gamma^n_2=\bar
w_{[-r_n,r_n]}\bigl(C_*,\bar w_u\bigl(P^{*,n},\Del\bigr)
\bigr),\\
 \gamma^n_3&=&2C_*(r_n)\Del,
\qquad
r_n=\bigl\|P^{*,n}\bigr\|_u+\bar c.
\end{eqnarray*}

Since $P^{*,n}$ converge uniformly, a.s., we have that $r_n$ converge
to a finite-valued r.v.
Since $C_*$ is a continuous function and $n^{-1/2}\varrho_n\to0$
\eqref
{51}, $\gamma^n_1\to0$ a.s.
as $n\to\iy$. Moreover, by the uniform convergence of $P^{*,n}$ to a
process with continuous sample paths,
we have $\lim_{\Del\to0}\limsup_{n\to\iy}\gamma^n_2=0$. It
follows that
for some $\gamma^n_4=\gamma^n_4(\Del)\ge0$
satisfying $\lim_{\Del\to0}\limsup_{n\to\iy}\gamma^n_4=0$,
%
%
\begin{eqnarray}
\label{73} \int_0^uC_*\bigl(E^n(t)
\bigr)\,dt&\ge&\int_0^uC_*\bigl(P^{*,n}(t)
\bigr)\,dt+\bigl[C_*(c\varrho _n)-K^n\bigr](\til
\tau_n-\tau_n)
\nonumber
\\[-8pt]
\\[-8pt]
\nonumber
&&{} -K^n1_{\{\tilde{\tau}^n<u\}}-
\gamma^n_4,
\end{eqnarray}
holds on $\Om^n$, where
\[
K^n=\sup_{t\in[0,u]}C_*\bigl(P^{*,n}(t)\bigr).
\]
Now, $\limsup_nK^n$ is a finite r.v.
On the other hand, one has $\liminf_{r\to\iy}C_*(r)/\break r>0$ due to the
fact that $C_*$
is convex, nondecreasing and nonconstant
whence (recalling that $\varrho_n\tinf$ as $n\tinf$) $C_*(c\varrho
_n)\to
\iy$ grows without bound. Thus the second
term on the RHS of \eqref{73}
is negative for only finitely many $n$. Also, if $\tilde{\tau}^n<u$,
then we have, by definition, that $\tilde{\tau}^n-\tau^n=t_n$ and
(recalling that $t_n\varrho_n\tinf$) we have that the second term on
the RHS of \eqref{73} grows to infinity so that the sum of the second
and third terms in \eqref{73} is negative for only finitely many $n$.
Since $1_{\Om^n}\to1$ a.s., we have thus shown that a.s.,
\[
\liminf_{n\to\iy}\int_0^uC_*
\bigl(E^n(t)\bigr)\,dt\ge\liminf_{n\to\iy}\int
_0^uC_*\bigl(P^{*,n}(t)\bigr)\,dt =\int
_0^uC_*\bigl(Q^*(t)\bigr)\,dt.
\]
Noting that, by definition of $C_*$ and $E^n$, $C(\wh Q^n(t))\ge C_*(E^n(t))$,
and recalling that $u$ is arbitrary completes the proof.
\end{pf*}

\subsection{\texorpdfstring{Proof of Proposition \protect\ref{prop2}}
{Proof of Proposition 4.1}}\label{sec42}

The result is proved in three major steps, where the first two
establish the lower bounds \eqref{70} and \eqref{71}
for suitably defined processes $P^{*,n}$ and $H^{*,n}$
[note that \eqref{72} is immediate from the nonnegativity of $\wh
Q^n$]. The third
step verifies the statements of Proposition \ref{prop2} with regard to
convergence.

\textit{Step \textup{1:} The interval $[0,\tau_n)$.}
In this step we analyze the time interval $[0,\tau_n)$, introduce
processes $P^{*,n}$ and $H^{*,n}$ [cf. \eqref{74} and \eqref{75}]
and argue that they satisfy the bound \eqref{70}.

We will next need a result that is similar to the minimality property
of the Skorohod map $\Gam$ \eqref{54} but that allows for a certain
kind of perturbation.
The minimality property of the Skorohod map is stated as follows: Let
$\zeta\in\mathcal{D}$. Let $\eta\in\mD$ be nondecreasing and satisfy
$\eta(0)\ge0$.
Assume $\zeta(t)+\eta(t)\ge0$, for all $t\ge0$. Then
%
%
\begin{equation}
\label{55} \zeta(t)+\eta(t)\ge\Gam[\zeta](t)\equiv\zeta(t)+\sup
_{s\le
t}\bigl[\zeta (s)^-\bigr],\qquad t\ge0.
\end{equation}

The following lemma provides the variant that we need.

\begin{lemma}
\label{lem3}
Let $u>0$ and $\eps>0$, $\eps<u$, be given.
Let $\zeta\in\mD$ and assume $\zeta(0)\ge0$. Let
\[
\al=\zeta+\eta+\beta,
\]
where $\eta\in\mD$ is nondecreasing and satisfies $\eta(0)\ge0$,
$\beta
\in\mD$ satisfies
%
%
\begin{equation}
\label{58} -\eps^2\le\int_0^t
\beta(s)\,ds\le\eps^2,\qquad t\in[0,u]
\end{equation}
and $\al(t)\ge0$, $t\in[0,u]$.
Then
%
%
\begin{equation}
\label{56} \al(t)\ge\Gam[\zeta](t)+\beta(t)-\bar w_u(\zeta,\eps)-3
\eps,\qquad t\in[0,u].
\end{equation}
\end{lemma}

\begin{pf}
By \eqref{55}, we have $\al\ge\Gam[\zeta+\beta]$. Thus
\[
\al(t) \ge\zeta(t)+\beta(t) + \sup_{s\le t}\bigl[\bigl(\zeta(s)+
\beta(s)\bigr)^-\bigr].
\]
Denote $\delta=\bar w_u(\zeta,\eps)+3\eps$.
To prove the claim, it suffices to show that
%
%
\begin{equation}
\label{62} \sup_{[0,t]}\bigl[(\zeta+\beta)^-\bigr]\ge\sup
_{[0,t]}\bigl[\zeta^-\bigr]-\delta,
\end{equation}
because then $\al\ge\zeta+\beta+\sup_{[0,\cdot]}[\zeta^-]-\delta
=\Gam
[\zeta]+\beta-\delta$.
To this end, consider first $t\in[0,\eps]$.
Since $\zeta(0)\ge0$ by assumption, we have
\[
\sup_{[0,\eps]}\bigl[\zeta^-\bigr]-\delta\le\bar w_u(
\zeta,\eps)-\delta =-3\eps.
\]
Thus \eqref{62} is immediate from the nonnegativity of the LHS of that
inequality.

Next, fix $t\in[\eps,u]$. Toward showing that \eqref{62} holds in this
case as well,
note that, for any $s\in[\eps,u]$,
%
%
\begin{equation}
\label{57} \inf_{[s-\eps,s]}\beta\le3\eps
\end{equation}
for otherwise we would have, by \eqref{58},
\[
3\eps^2\le\int_{s-\eps}^s\beta(\tau)\,d
\tau\le2\eps^2.
\]
(note that \eqref{58} is imposed for all $t\in[0,u]$ and, in
particular, for $t=s-\varepsilon$---we are using that here). We
consequently have
\[
\inf_{[s-\eps,s]}(\zeta+\beta)\le\sup_{[s-\eps,s]}\zeta+
\inf_{[s-\eps,s]} \beta\le\sup_{[s-\eps,s]}\zeta+3\eps \le
\inf_{[s-\eps,s]}\zeta+\delta,
\]
where in the second inequality we used \eqref{57} and our choice of
$\delta$ above.
Taking the infimum over $s\in[\eps,t]$, we obtain
%
%
\begin{equation}
\label{59} \inf_{[0,t]}(\zeta+\beta)\le\inf_{[0,t]}
\zeta+\delta.
\end{equation}
Let us deduce from the above that
%
%
\begin{equation}
\label{60} \inf_{[0,t]}0\w(\zeta+\beta)\le\inf
_{[0,t]}0\w\zeta+\delta.
\end{equation}
Indeed, if $\inf_{[0,t]}\zeta\ge-\delta$, then the RHS of \eqref
{60} is
nonnegative, and hence
this inequality is valid.
If $\inf_{[0,t]}\zeta<-\delta$, then $\inf\zeta=\inf0\w\zeta$,
and hence
the claim follows \eqref{59}. In both cases,
\eqref{60} holds. Note that \eqref{60} is equivalent to \eqref{62}. We
have thus shown that \eqref{62} holds
for any $t\in[0,u]$.
This completes the proof.
\end{pf}

We proceed with the proof of Theorem \ref{th1}. Recall equation \eqref
{31} for $\wh Q^n$ and that $E^n=\theta'\wh Q^n$. Writing
%
%
\begin{eqnarray}
\label{61} F^n&=&\theta'\wh X^n(0)+
\theta'W^n-\int_0^\cdot
\sum_{i,j}\theta _i\eps
^n_{ij}\til B^n_{ij}(s)\,ds,
\\
\label{64} G^n&=&-\int_0^\cdot\sum
_{i,j}\theta_i\mu_{ij}\til
B^n_{ij}(s)\,ds, \qquad H^n=-\sum
_{i,j}\theta_i\wh B^n_{ij},
\end{eqnarray}
we have
%
%
\begin{equation}
\label{63} E^n=F^n+G^n+H^n.
\end{equation}
We will apply Lemma \ref{lem3}, substituting $E^n(t)$, $F^n(t)$,
$G^n(t)$ and $H^n(t)$, $t\in[0,\tau_n]$,
for $\al$, $\zeta$, $\eta$ and $\beta$, respectively.
To this end, let us verify the assumptions on these processes.
First, by \eqref{20}, \eqref{21} and \eqref{18}, we have $W^n(0)=0$.
Since $\wh X^n(0)$ is assumed to have values in
$\R_+^\bI$, we have by \eqref{61} that $F^n(0)\ge0$.

Let us show that $G^n$ is nondecreasing. It suffices to show that, for
fixed $n$ and $t$,
%
%
\begin{equation}
\label{66} \sum_{i,j}\theta_i
\mu_{ij}\til B^n_{ij}(t)\le0.
\end{equation}
We do this by invoking Lemma \ref{lem1}. Let $x\in\Sig$ be defined by
$x_{ij}=-\nu_j^{-1}\til B^n_{ij}(t)$.
By~\eqref{03} and the fact that $\xi^*_{ij}=0$ for $(i,j)\in\calE
_{nb}$, we have $x_{ij}\le0$
for $(i,j)\in\calE_{nb}$. Also, for any $j$, $\sum_ix_{ij}=-\nu
_j^{-1}\sum_i\til B^n_{ij}(t)=I^n_j(t)$, by \eqref{12}.
Hence $\sum_ix_{ij}\ge0$. Thus by Lemma \ref{lem1}, $\theta'\bar
\mu
(x)\ge0$. Recalling that
$\bar\mu_{ij}=\nu_j\mu_{ij}$, we obtain \eqref{66}.

Clearly, $G^n(0)=0$.

Let us show that $H^n$ satisfies a bound of the form \eqref{58}.
By the definition of $\tau_n$, $|\int_0^t\til B^n_{ij}(s)\,ds|\le
\varrho
_n$, for all $t\le\tau_n$.
Hence by \eqref{03},
\[
\biggl|\int_0^t\wh B^n_{ij}(s)\,ds
\biggr|\le n^{-1/2}\varrho_n, \qquad t\in [0,\tau_n]
\]
and
%
%
\begin{equation}
\label{69} \biggl|\int_0^t H^n(s)\,ds \biggr|\le
c_0n^{-1/2}\varrho_n,\qquad t\in [0,
\tau_n],
\end{equation}
some constant $c_0$.

Finally, note that $E^n(t)\ge0$ for all $t$, since $\wh Q^n_i(t)\ge0$
for every $i$ and all $t$.

Having verified the hypotheses of Lemma \ref{lem3}, we obtain
$E^n(t)\ge\Gam[F^n](t)-\delta_n+H^n(t)$, for all $t\in[0,\tau_n]$, where
%
%
\begin{equation}
\label{77} \delta_n=\bar w_{\tau_n}(F_n,\bar
\eps_n)+3\bar\eps_n,\qquad \bar\eps_n=c_0^{1/2}n^{-1/4}
\varrho_n^{1/2}.
\end{equation}
Set
%
%
\begin{eqnarray}
\label{74} P^{*,n}(t)&=&\bigl(\Gam\bigl[F^n\bigr](t)-
\delta_n\bigr)1_{\{t<\tau_n\}}+\Gam\bigl[\theta'\wh
X^n(0)+\theta'\wh W^n\bigr](t)1_{\{t\ge\tau_n\}},
\\
\label{75} H^{*,n}(t)&=&H^n(t)1_{\{t<\tau_n\}}.
\end{eqnarray}
Then $(P^{*,n},H^{*,n})$ agree with $(\Gam[F^n]-\delta_n,H^n)$
on the interval $[0,\tau_n)$, and we have shown
%
%
\begin{equation}
\label{65} E^n(t)\ge P^{*,n}(t)+H^{*n}(t), \qquad t
\in[0,\tau_n).
\end{equation}

\textit{Step \textup{2:} The interval $[\tau_n,\til\tau_n)$.} In this step we
show that \eqref{71} holds
on a suitably defined event $\Om^n$.
The argument is based on the following lemma.
For $x\in\Sig$ denote $\|x\|^2=\sum_{i,j}x_{ij}^2$. For the lemma below
recall that
$\Sig$ is the set of $\bI\times\bJ$ matrices for which the $(i,j)$
entry is zero whenever $i\nsim j$. Also, note that in this lemma the
uniqueness of $\xi^*$ is used in a crucial manner.
%
\begin{lemma}
\label{lem4}
Let $x\in\Sig$ be such that $x_{ij}\le0$ for all $(i,j)\in\calE
_{nb}$, and
$\sum_ix_{ij}\ge0$ for all $j\in\calJ$.
Then
\[
\max_i\bar\mu_i(x)\ge c_1\|x\|,
\]
where $c_1>0$ is a constant that does not depend on $x$.
\end{lemma}
\begin{pf}
Let
\[
K=\biggl\{\xi^*-x\dvtx\|x\|=\eps, x_{ij}\le0, (i,j)\in
\calE_{nb}, \sum_ix_{ij}\ge0,j
\in\calJ\biggr\}.
\]
A review of the proof of Lemma \ref{lem1} shows that
$K\subset\X$, provided that $\eps>0$ is sufficiently small.
Let such $\eps$ be fixed. Recall that $\lambda=\bar\mu(\xi^*)$,
and note
that $\xi^*\notin K$.
Thus the uniqueness of $\xi^*$, stated in Assumption \ref{assn2},
implies that
there is no $\xi\in K$ for which $\bar\mu(\xi)=\lambda$. Hence
$\lambda\notin
\bar\scrM$, where $\bar\scrM$
is the image $\bar\mu(K)$ of $K$ under $\bar\mu$.
Recall that $\lambda$ is a maximal element in $\scrM$ with respect to the
usual partial order
in~$\R^\bI$. Since $\lambda\notin\bar\scrM$, this says
\[
\max_i\bigl[\lambda_i-\bigl(\bar\mu(\xi)
\bigr)_i\bigr]>0\qquad \mbox{for every } \xi\in K.
\]
Note that $K$ is a compact set, and that the LHS of the above display
depends continuously on $\xi$. Hence there exists $\delta>0$ such that
\[
\max_i\bigl[\lambda_i-\bigl(\bar\mu(\xi)
\bigr)_i\bigr]\ge\delta \qquad\mbox{for every } \xi\in K.
\]
Noting that the conclusion of the lemma holds for $x=0$,
consider any nonzero member $x$ of $\Sig$, satisfying the lemma's assumptions.
Then $\xi^*-\eps\|x\|^{-1}x\in K$.
Hence
\begin{eqnarray*}
\eps\|x\|^{-1}\max_i\bar\mu_i(x) &=&
\max_i\bar\mu_i\bigl(\eps\|x\|^{-1}x
\bigr)
\\
&=&\max_i\bigl[\bar\mu\bigl(\xi^*\bigr)-\bar\mu\bigl(\xi^*-
\eps\|x\|^{-1}x\bigr)\bigr]
\\
&=&\max_i\bigl[\lambda_i-\bar\mu\bigl(\xi^*-
\eps\|x\|^{-1}x\bigr)\bigr]
\\
&\ge&\delta.
\end{eqnarray*}
The claim follows with $c_1=\delta/\eps$.
\end{pf}

We have already argued that if $x\in\Sig$ is defined by setting
$x_{ij}=-\nu_j^{-1}\til B^n_{ij}(t)$,
then $x$ satisfies the assumptions of Lemma \ref{lem1}, equivalently
Lemma \ref{lem4}. As a consequence,
so does $y:=\int_0^{\tau_n}x(t)\,dt$. By \eqref{32}, on the event
$\tau
_n<u$, we have that, for some $(i,j)\in\calE$,
$|y_{ij}|=\nu_j^{-1}\varrho_n\ge\varrho_n$, and so $\|y\|\ge
\varrho_n$.
Applying Lemma \ref{lem4} yields that
there exists $i^*\in\calI$ such that
\[
\bar\mu_{i^*}(y)\ge c_1\|y\|\ge c_1
\varrho_n.
\]
Namely,
\[
\sum_j\mu_{i^*j}\int
_0^{\tau_n}\til B^n_{i^*j}(t)\,dt
\le-c_1\varrho_n.
\]
Invoking \eqref{31} and using the nonnegativity of $\wh Q^n_i(0)$, we obtain
%
%
\begin{equation}
\label{67} \wh Q^n_{i^*}(\tau_n)\ge
W^n_{i^*}(\tau_n)+c_1
\varrho_n-\sum_j\eps^n_{i^*j}
\varrho_n-\sum_j\wh
B^n_{i^*j}(\tau_n),
\end{equation}
where we used the fact that $|\int_0^{\tau_n}\til B^n_{ij}(t)\,dt|\le
\varrho_n$, by \eqref{32}.
By \eqref{01}, $\sum_jN^n_j\le c_2n^{1/2}$, for some constant
$c_2>0$. Hence
%
%
\begin{eqnarray}
\label{76}\qquad \sum_j\wh B^n_{ij}(t)
\le n^{-1/2}\sum_j\til
B^n_{ij}(t)\le n^{-1/2}\sum
_j B^n_{ij}(t)\le n^{-1/2}
\sum_jN^n_j\le
c_2,
\nonumber
\\[-8pt]
\\[-8pt]
\eqntext{i\in\calI.}
\end{eqnarray}
Also, by \eqref{51}, for $n\ge n_0$, the third term on the RHS of
\eqref
{67} is bounded by 1,
some deterministic $n_0$.
As a result, for a suitable constant $c_3>0$, we have, for $n\ge n_0$,
\[
\wh Q^n_{i^*}(\tau_n)\ge W^n_{i^*}(
\tau_n)+c_3\varrho_n.
\]
Recall $\eps^n_M=\max_{i,j}|\eps^n_{ij}|$ and let $t_n=c_3\rho
_nc_2^{-1}(\eps^n_M)^{-1}n^{-1/2}/2$
(where $t_n=\iy$ if $\eps^n_M=0$). Using \eqref{31}, and the bound
\[
\sum_j\eps^n_{ij}\til
B^n_{ij}\le c_2\eps^n_Mn^{1/2},
\]
we have
\begin{eqnarray*}
\wh Q^n_{i^*}(t) &\ge& W^n_{i^*}(
\tau_n)+c_3\varrho_n-c_2
\eps^n_Mn^{1/2}(t-\tau_n)
\\
&\ge& W^n_{i^*}(\tau_n)+\frac{c_3}{2}
\varrho_n,\qquad t\in[\tau_n, \til\tau_n],
\end{eqnarray*}
where, as in the statement of Proposition \ref{prop2}, $\til\tau
_n=(\tau
_n+t_n) \w u$.
By the nonnegativity of $\wh Q^n_i$ and positivity of $\theta_i$,
$i\in
\calI$, we have, with $\theta_m=\min_i\theta_i$,
%
%
\begin{equation}
\label{68} E^n(t)\ge\theta_{i^*}\wh Q^n(t)\ge
\theta_m \biggl(\frac
{c_3}{2}\varrho _n-
\bigl\|W^n(\tau_n)\bigr\| \biggr),\qquad t\in[\tau_n, \til
\tau_n].
\end{equation}
This shows that \eqref{71} is valid on $\Om^n:=\{\|W^n(\tau_n)\|
<c_3\varrho_n/4\}\cap\{n\ge n_0\}$.

\textit{Step \textup{3:} Convergence}. We are now in a position to prove all
statements of the proposition.
By \eqref{69}, \eqref{70-} holds up to $\tau^n$; by \eqref{75},
this estimate remains valid up to $u$. Moreover,
it is shown in \eqref{76} that $\wh B^n$ are bounded above, and it can
similarly be shown that they
are bounded below. Hence $H^{*,n}(t)$ are bounded uniformly in $n$ and $t$.
Note that $\varrho_nt_n=c\varrho_n^2(\eps_M^n)^{-1}n^{-1/2}$, for some
constant $c$. By \eqref{51},
this product converges to $\iy$.
Thus to complete the proof, it remains to argue that $1_{\Om^n}\to1$
a.s. and $P^{*,n}\to Q^*$ a.s.

Toward this end, denote
%
%
\begin{equation}
\label{49} \bar T_{ij}(t)=\bar\mu_{ij}\xi^*_{ij}t,\qquad
t\ge0, (i,j)\in \calE.
\end{equation}
Let us first show that
%
%
\begin{equation}
\label{50} \bigl|\bar T^n_{ij}-\bar T_{ij}\bigr|_{\tau_n}
\to0\qquad \mbox{in probability}.
\end{equation}
Indeed, by \eqref{45}, \eqref{03} and \eqref{35},
\[
\bar T^n_{ij}(t)=n^{-1}\mu^n_{ij}
\xi^*_{ij}N^n_jt+n^{-1}\mu
^n_{ij}\int_0^t\til
B^n_{ij}(s)\,ds =\bar\mu^n_{ij}
\xi^*_{ij}t+n^{-1}\mu^n_{ij}\int
_0^t\til B^n_{ij}(s)\,ds.
\]
Hence by \eqref{32} and \eqref{49},
%
%
\begin{equation}
\label{52} \bigl|\bar T^n_{ij}(t)-\bar T_{ij}(t)\bigr|\le
\wh\delta_n:=\bigl|\bar\mu ^n_{ij}-\bar\mu
_{ij}\bigr|\xi^*_{ij}u+n^{-1}\mu^n_{ij}
\varrho_n,\qquad  t\in[0,\tau_n].
\end{equation}
Using the convergence \eqref{35} and
recalling that $\mu^n_{ij}$ is asymptotic to $\mu_{ij}n^{1/2}$, it
follows by \eqref{51} that $\wh\delta_n$, defined in the above
display, converges to zero. This shows \eqref{50}.

Recall the process $W$ from Lemma \ref{lem2}(iii). We now argue that
%
%
\begin{equation}
\label{53} \bigl\|W^n-W\bigr\|_{\tau_n}\to0 \qquad\mbox{a.s.}
\end{equation}
By \eqref{21}, \eqref{18} and \eqref{46},
\begin{eqnarray*}
W^n_i(t)-\wh W^n_i(t)&=&-\sum
_j\bigl[\wh S^n_{ij}\bigl(
\bar T^n_{ij}(t)\bigr)-\wh S^n_{ij}
\bigl(\bar\mu_{ij}\xi^*_{ij}t\bigr)\bigr] \\
&=&-\sum
_j\bigl[\wh S^n_{ij}\bigl(\bar
T^n_{ij}(t)\bigr)-\wh S^n_{ij}\bigl(
\bar T_{ij}(t)\bigr)\bigr].
\end{eqnarray*}
Hence by \eqref{52},
\[
\bigl|W^n_i(t)-\wh W^n_i(t)\bigr|\le\sum
_j\bar w_{u_1}\bigl(\wh
S^n_{ij},\wh\delta _n\bigr),\qquad t\in[0,
\tau_n],
\]
where $u_1=\max_{(i,j)\in\calE}\bar T_{ij}(u)+1$. By Lemma~\ref
{lem2}(i), the processes $\wh S^n_{ij}$ converge,
uniformly on compacts, to processes with continuous sample paths. Hence
the RHS of the above display
converges to zero. Applying Lemma \ref{lem2}(iii), we conclude \eqref{53}.

It follows from \eqref{53} that $1_{\Om^n}\to1$ a.s.

Next, by definition of $\Gam$, for any $v>0$, the mapping
$x|_{[0,v]}\mapsto\Gam[x]|_{[0,v]}$ is Lipschitz
continuous in the sup norm, with constant $2$. Recalling the definition
of $P^{*,n}$ \eqref{74},
we obtain, for some constant $c$,
\[
\bigl\|P^{*,n}-Q^*\bigr\|_u\le|\delta_n|+c
\bigl\|F^n-\theta'X_0-\theta'W
\bigr\|_{\tau
_n} +c\bigl\| \wh X^n(0)-X_0\bigr\|+c\bigl\|\wh
W^n-W\bigr\|_u.
\]
We have already argued that the last two terms above converge to zero
a.s. [Lemma~\ref{lem2} and Remark \ref{rem1}(a)].
By \eqref{61}, for some constant $c_1$,
\[
\bigl\|F^n-\theta'X_0-\theta'W
\bigr\|_{\tau_n}\le c_1\bigl\|\wh X^n(0)-X_0
\bigr\|+c_1\bigl\| W^n-W\bigr\|_u +c_1
\eps^n_M\varrho_n\to0\qquad \mbox{a.s.},
\]
where we used \eqref{51} and \eqref{53}.
By \eqref{51} and \eqref{77}, $\bar\eps_n\to0$, hence the convergence
in the above display
implies that $\delta_n\to0$ a.s. This shows the convergence of $P^{n,*}$
to~$Q^*$,
and completes the proof.

\section{Asymptotic optimality in heavy-traffic}\label{sec:asymopt}

\subsection{The tracking policy}

In this section we devise a sequence of controls that asymptotically
achieve the lower bound in Theorem \ref{th1}. To construct them and
prove their asymptotic optimality we need some further assumptions.
Recall from~\eqref{eq:cstartdefin} that
%
%
\begin{equation}
\label{eq:2} C_*(a)=\inf\bigl\{C(q)\dvtx q\in\R_+^\bI,
\theta'q=a\bigr\}.
\end{equation}

\begin{assumption}[(Continuous minimizer)]
There exists a
locally Lipschitz function $f\dvtx\R_+ \to\R_+^I$ such that $\theta
'f(a)=a$ and $C(f(a))=C_*(a)$ for all $a\geq0$. \label
{asum:continuousselection}
\end{assumption}
We extend the function $f$ to the real line by setting it to zero on
$(-\iy,0)$ (note
that the extended function is continuous).
In fact, the actual implementation of the policy will be based on small
perturbations of $f$, that we denote by $f^n$. These perturbations are
explicitly provided below; see \eqref{eq:qndefin}.

The policies we construct will not use the nonbasic activities at all.
It will be
convenient, in terms of notation, to disregard these activities by
assuming they do not
exist. Thus we assume (w.l.o.g., as far as the results of this section
are concerned) that
\textit{all activities in the model are basic}. In particular, $\calG
=\calG
_b$, and $i\sim j$ means $(i,j)$
is an activity, equivalently, a basic activity. Since Assumptions~\ref
{assn2} and~\ref{assn3} are still in force, the graph is a tree.

Based on the structure of the control which achieves the lower bound
for the \textit{diffusion} control problem
(Proposition \ref{prop1}), we seek
control policies for the queueing model having two main properties.
Namely, (i) that
the sequence $\theta' \hQ^n$ converges in law to the RBM $Q^*$,
and (ii) that, given $\theta'\wh Q^n=a$,
$\wh Q^n$ itself is close to the minimizing $f$ in \eqref{eq:2}.
To be more precise about (ii), denote
%
%
\begin{equation}
\label{86} \cX^n=f^n\bigl(\theta'
\hX^n\bigr) \quad\mbox{and}\quad \cQ^n=f^n\bigl(
\theta'\hQ^n\bigr).
\end{equation}
Further, note by \eqref{13} that $\|\wh X^n-\wh Q^n\|\le\|\wh B^n\|$.
Then property (ii) corresponds to having $\wh Q^n-\check Q^n\To0$.
If $\wh B^n\To0$, then the above can be achieved by having
$\wh Q^n-\check X^n\To0$.
It turns out to be more convenient to work with the latter,
that is, to prove $\wh B^n\To0$ and $\wh Q^n-\check X^n\To0$.

The proposed policy, to which we refer as \textit{the tracking policy}, seeks
to achieve the convergence $\wh Q^n-\check X^n\To0$
by letting $\wh Q^n$ track $\check X^n$. That is,
upon service completion in pool $j$ at time $t$, the policy assigns to
the newly-available
server a customer from a class within the set
%
%
\begin{equation}
\bigl\{i\dvtx i\sim j, \hQ_i^n(t-)>\check
X^n_i(t-) \bigr\}, \label{eq:interim0}
\end{equation}
so as to decrease the difference.
This, however, is not a precise description of the policy, as
the choice of the class for service in \eqref{eq:interim0} will not be
arbitrary.
It will rely on the structure of the tree. For a precise statement
we need some additional notation.

Let $\calV=\calI\cup\calJ$ denote the vertex set of $\calG$. Pick any
$i_0\in\calI$ and designate it as the root. Denote by $d(k)$ the graph
distance of a node $k\in\calV$ from the root.
For $i\in\calI\setminus\{i_0\}$ denote by $\jbar(i)$ the neighbor
$j\sim i$ that is closer to the root
than $i$, and by $\calJ(i)$ the (possibly empty) set of neighbors
$j\sim i$ that are farther
from the root than $i$.
We sometimes refer to $\jbar(i)$ and $\calJ(i)$ as the nodes right
\textit{above} and, respectively, \textit{below} $i$, thinking of the tree as
being depicted with
the root at the top.
For nodes $j\in\calJ$, define analogously $\ibar(j)$ and
$\calI(j)$. Next, fix a labeling of all graph nodes by distinct numbers
between $1$ and $\bI+\bJ$,
so that members in $\calI$ have the labels $\{1,2,\ldots,\bI\}$ and
those in $\calJ$
have the labels $\{\bI+1,\bI+2,\ldots,\bI+\bJ\}$.
Identify the set $\calV$ with $\{1,\ldots,\bI+\bJ\}$ accordingly.
The labeling should satisfy the following additional condition, namely
for $i_1,i_2\in\calI$:
\[
d(i_1)<d(i_2) \mbox{ implies } i_1>i_2,
\]
and for $j_1,j_2\in\calJ$,
\[
d(j_1)<d(j_2) \mbox{ implies } j_1>j_2.
\]
We let $j_0=\max\{j\dvtx j\sim i_0\}$. Some of our statements preclude the
root $i_0$. We write $\I_{-i_0}=\I\setminus\{i_0\}$ and, for a vector
$x\in\mathbb{R}^{\bI}$, we write $x_{-i}=(x_1,\ldots,x_{i-1},\break x_{i+1},\ldots, x_{\bI})$.

\textit{Perturbed functions}: As alluded to earlier, we will work with a
perturbed version, $f^n$, of the function $f$. Fix a function $f$ that
satisfies Assumption \ref{asum:continuousselection}. Let sequences $\{
\kappa_n\}$ and $\{\bar\kappa_n\}$ with $\kappa_n/\bar\kappa
_n\rightarrow0$ and $\bar\kappa_n\rightarrow0$ be given. For $i\in
\I
_{-i_0}$ we set
%
%
\begin{equation}
f_i^n(x)=\cases{ %
( \bI
\theta_i)^{-1}x,&\quad $x\in[0,\kappa_n),$
\vspace*{2pt}\cr
(\bI\theta_i)^{-1}\kappa_n,&\quad $x\in[
\kappa_n,\bar\kappa_n),$
\vspace*{2pt}\cr
f_i(x) (1-\bar\kappa_n/x)+(\bI\theta_i)^{-1}
\kappa_n,&\quad $x\in[\bar \kappa _n,\infty).$}
\label{eq:qndefin}
\end{equation}
Also set $f_{i_0}^n(x)=(x-\sum_{i\neq
i_0}\theta_if_i^n(x))/\theta_{i_0}$.

\begin{remark}[(Properties of $f^n$)]
The perturbed
functions $f^n$ are such that, for $i\in\I_{-i_0}$, $f_i^n$ is small
but strictly positive in the vicinity of $0$. Thus,
if one can guarantee that $\hQ_i^n\approx f_i^n(\theta'\hX^n)$, then
$\hQ_i^n>0$ for all $i\in\I$ whenever $\theta'\hX^n>0$, thus there are
no idling servers. This property is used
in the proof.
Additional observations regarding $f^n$ will be useful in what
follows: $f_i^n(x)\geq\kappa_n$ for $x\geq\kappa_n$ and $i\in\I
_{-i_0}$. Also, $f_{i_0}^n(x)\geq(\bI\theta_{i_0})^{-1}x$ for all
$x\in[0,\bar\kappa_n)$ and $f_{i_0}^n(x)\geq(\theta
_{i_0})^{-1}\bar
\kappa_n\geq(\bI\theta_{i_0})^{-1}\bar\kappa_n$ for all $x\geq
\bar
\kappa_n$. Since $\theta_if_i^n(x)\leq x$ for all $i\in\I$ and
$x\geq
0$, we have that
%
%
\begin{eqnarray}
\qquad\sup_{x\geq0}\bigl|\theta_if_i^n(x)-
\theta_if_i(x)\bigr|\leq\bar \kappa _n+\sup
_{x\geq\bar\kappa_n}\bigl|\theta_if_i^n(x)-
\theta_if_i(x)\bigr|\leq 2\bar \kappa_n
\rightarrow0,
\nonumber
\\[-8pt]
\\[-8pt]
 \eqntext{i\in\I.}\label{eq:qconv}
\end{eqnarray}
Finally, it is easy
to verify that $f^n$ are locally Lipschitz uniformly in $n$.
\label{rem:qnproperties}
\end{remark}

Recall that $I_j^n(t)$ is the number of idle servers in pool $j$ at
time $t$.

\textit{The tracking policy}:
\begin{longlist}[(ii)]
\item[(i)]
Upon each arrival of a class-$i$ customer, say at time $t$,
if there are idle servers in one of the pools
$j\in\calJ(i)$, then it is routed to the pool
\[
\min\bigl\{j\in\calJ(i)\dvtx I_j^n(t-)>0\bigr\}.
\]
Otherwise it is queued [even if there are idle servers in pool $\jbar(i)$].
\item[(ii)]
Upon each service completion in pool $j$, say at time $t$,
the server admits to service a customer from class
\[
\min\mathcal{K}_j(t-),
\]
where
\[
\mathcal{K}_j(t-):= \bigl\{k\in\I(j)\dvtx\hQ_k^n(t-)>
\cX _k^n(t-) \bigr\},
\]
provided this set is nonempty.
If $\mathcal{K}_j(t-)=\varnothing$ but $\hQ_{\ibar(j)}^n(t-)>0$ the
server admits to service a customer from class $\ibar(j)$. Otherwise,
the server remains idle.
\end{longlist}
%

\begin{remark}[(Work conservation)]\label{rem:workconservation}
It is clear that the policy is not work conserving.
However, it is not hard to see that
at any given time $t$,
all servers in pool $j$ must be busy at $t$, provided $Q_i^n(t)>0$ for
\textit{all} $i\sim j$.
Note that this means $\sum_{i}B_{ij}^n(t)=N_j$ [and, in turn,
$\sum_{i}\tilde{B}_{ij}^n(t)=0$]. This property will be useful in
what follows.
\end{remark}

The main result of this section is the following.
For simplicity, the initial conditions for $\wh B^n$ and
$\wh Q^n$ are assumed to vanish.
%
\begin{theorem} Suppose that Assumptions \ref{assn2}, \ref{assn3} and
\ref{asum:continuousselection} hold and that $\hB^n(0)=\hQ^n(0)=0$ for
all $n$. Then, under the tracking policy,
%
%
\begin{equation}
\bigl(\hQ^n-\cQ^n,\hB^n,
\theta'\hQ^n\bigr)\Rightarrow\bigl(0,0,Q^*\bigr),
\label{eq:conv_track}
\end{equation}
where $Q^*$ is as in Theorem \ref{th1}. Consequently,
\[
\int_0^uC\bigl(\wh Q^n(t)\bigr)\,dt
\Rightarrow\int_0^u C_*\bigl(Q^*(t)\bigr)\,dt,
\]
where $C_*(\cdot)$ is as in \emph{(\ref{eq:cstartdefin})}. \label
{thm:conv_track}
\end{theorem}

\begin{remark}[(Special cost structures)]
\textit{Separable convex costs}:
Consider, in terms of costs, the setting of \cite{mansto}. This is the
case that $C(q)=\sum_{i\in\I}C_i(q_i)$ where $C_i, i\in\I$, are twice
continuously differentiable strictly increasing and strictly convex
functions with $C_i'(0)=0$ for all $i\in\I$. Fixing the constant $a$,
under the Kuhn--Tucker conditions for \eqref{eq:2},
the unique solution $f(a)$ must satisfy
\[
C_i'\bigl(f_i(a)\bigr)=-y(a)
\theta_i,
\]
where $y(a)$ is the Lagrange multiplier of the constraint $\theta'q=a$.
Thus
$f_i(a)=(C_i')^{-1}(-y(a)\theta_i)$. It can be verified that $f(a)$ is
locally Lipschitz
and, in particular, that it satisfies Assumption \ref
{asum:continuousselection}. Moreover, since $C_i$ is strictly
increasing and strictly convex, $f_i(a)$ is a strictly increasing function.

Recall from
\eqref{eq:zdefin} that $\theta_i\mu_{ij}=z_j^*$ for all $(i,j)\in
\calE
_b$. The
Kuhn--Tucker
condition is equivalently written as
\[
\mu_{ij}C_i'\bigl(f_i(a)
\bigr)=-y(a)z_j^*\qquad\mbox{for all } (i,j)\in\calE_b,\vadjust{\goodbreak}
\]
and, in particular, for all $a\geq0$,
\[
\mu_{ij}C_i'\bigl(f_i(a)
\bigr)=\mu_{kj}C_k'\bigl(f_k(a)
\bigr)\qquad\mbox{for all } j\mbox{ and } i,k\in\I(j)\cup\bar{i}(j).
\]
Thus, our policy is consistent with the $Gc\mu$ policy in \cite{mansto}
in that it
aims at setting
\[
\mu_{ij}C_i'\bigl(\hQ^n_i(t)
\bigr)\approx\mu_{kj}C_k'\bigl(
\hQ_k^n(t)\bigr)\qquad\mbox{for all } j\mbox{ and } i,k\in
\I(j)\cup\bar{i}(j).
\]
The actual implementation is, however, different. Our service mechanism
follows the tree structure which is contrasted with the
$Gc\mu$ rule that would serve upon service completion a class in the set
\[
\mathop{\argmax}_{i\in\I(j)\cup\bar{i}(j)} \mu_{ij}C_i'\bigl(
\hQ^n_i(t)\bigr).
\]

\textit{Linear costs}: Suppose that $C(q)=\sum_{i}c_iq_i$ where
$c_i, i\in
\I$, are positive constants with $c_1\geq c_2\geq\cdots\geq c_{\bI}$, and
let us designate class $\bI$ as the root of the tree. In this
case, \eqref{eq:2} has a trivial solution $f_i(x)=0$ for all $i<{\bI}$
and $f_{\bI}(x)=x/\theta_{\bI}$.
Since $\kappa_n\rightarrow0$ as $n\tinf$, Theorem \ref{thm:conv_track}
guarantees that
$\hQ_i^n\Rightarrow0$ for all $i<\bI$ and $\theta_{\bI}\hQ_{\bI
}^n-\theta'\hQ^n\Rightarrow0, $ so that all queues except for the
lowest-cost queue are close to zero at diffusion
scale.

Our tracking policy is here a tree-based threshold policy as is the one
studied in \cite{belwil2}. If $t$ is such that $\theta'\hX
^n(t)>\kappa
_n$, then available servers give priority to classes $i<\bI$ that
exceed their threshold, that is, with $\hQ_i^n(t)\geq(\bI\theta
_i)^{-1}\kappa_n$. 
\label{rem:cost_examples}
\end{remark}

In the rest of this section we prove Theorem \ref{thm:conv_track}. All
symbols with superscript~$n$, such as $\wh Q^n$ and $\wh X^n$,
denote the respective processes under the tracking policy.

For sequences $a_n$ and $b_n$ of positive numbers, satisfying
$b_n/a_n>n^c$ for some $c>0$
and all large $n$, we write $a_n\ll b_n$.
Set $\varrho^n=n^{3/16}$ (satisfying the conditions we put on $\varrho
^n$ in the previous section). Recall that
$\varepsilon_M^n=\max_{i,j}|\varepsilon_{ij}^n|=O(n^{-1/4})$ (Lemma \ref
{lem2}), by which
$\varepsilon^n_M\varrho^n\ll1$.

By assumption, the inter-arrival times of all the processes
$(A_i, i\in\I)$ have finite moments of order $r>2$. Let $\alpha_A=(1/2-1/r)$.
Fix sequences $p_n$, $q_n$, $r_n$, $s_n$ of positive numbers, satisfying
%
%
\begin{equation}
\label{eq:gndefin} \eps^n_M\varrho^n\vee
n^{-\alpha_A}\ll p_n \ll q_n\ll r_n\ll
s_n\ll1.
\end{equation}

Assume, without loss of generality, that $p_n$ is given by $n^{-\alpha
_g}$ where $\alpha_g$
is a constant. Let $\omega\dvtx[0,\iy)\to[0,\iy)$ be given by $\omega
(x)=x^{\alpha
_{\omega}}$,
where $\alpha_\omega\in(1/3,1/2)$ is a constant. We further assume,
without loss of generality, that $\alpha_\omega>2\alpha_g$.

In \eqref{eq:qndefin} $\kappa_n$ and $\bar\kappa_n$ are such that
\[
p_n\ll\kappa_n\ll q_n\quad\mbox{and}\quad s_n\ll\bar\kappa_n\ll1.
\]

For the remainder of this section we fix the time horizon $u$. Recall
the processes
$\wh{A}^n$, $V^n$ and $W^n$, defined in \eqref{20}, \eqref{21} and
\eqref{18}, respectively.

\begin{lemma}
For every $\varepsilon>0$ there exist $K,L>0$ and $n_0$ such that
for $n\geq n_0$,
\[
\Pd \bigl\{\mbox{there exist } 0\leq s\leq t\leq u \mbox{ such that }\bigl \|
\La^n(t)-\La^n(s)\bigr\|>p_n+L\omega(t-s) \bigr\}
\leq \varepsilon
\]
and
\[
\Pd\bigl\{\bigl\|\La^n\bigr\|_u> K\bigr\}\le\varepsilon,
\]
where $\La^n$ is any one of the processes $\wh A^n$, $V^n$ and $W^n$.
\label{lem:Wtightness}
\end{lemma}

\begin{pf}
For the processes $\hat A^n$ and $\hat S^n$, the result follows
from strong approximations for renewal processes (see, e.g., Theorem
2.1.2 in \cite{csrgo1993weighted}) and the H\"{o}lder continuity of
Brownian motion paths.
For $V^n$ (and consequently for $W^n$), the result thus follows from
\eqref{21}, using the
uniform Lipschitz property of the processes~$\bar T^n_{ij}$
(note that $\mu_{ij}^nB_{ij}^n\leq cn$ where $c$ is constant).
\end{pf}

When fixing $\varepsilon>0$ we will, for simplicity of presentation,
assume that $n\geq n_0(\varepsilon)$.

Given $\varepsilon>0$, let $L=L_\eps$ be as in Lemma \ref{lem:Wtightness}.
Recall the process $G^n$ \eqref{64}.
We define the following random times:
\[
\sigma^n:=\inf \bigl\{t\geq0\dvtx\bigl\|\hB^n(t)\bigr\|\geq
s_n \bigr\}\wedge u
\]
and
\begin{eqnarray*}
&&\zeta^n:=\inf \bigl\{t\geq0\dvtx G^n(s_2)-G^n(s_1)
\geq r_n+ 4L\omega(s_2-s_1) \\
&&\hspace*{135pt}\mbox{for some }
0\leq s_1\leq s_2\leq t \bigr\} \wedge u.
\end{eqnarray*}

Finally, let $\tau^n$ be as in \eqref{32} and define $T^n=T^n(\eps)$
by
%
%
\begin{equation}
T^n:=\sigma^n\wedge\tau^n\wedge
\zeta^n. \label{eq:checktaudefin}
\end{equation}

\begin{proposition}\label{thm:conv_stopped}
Under the assumptions of Theorem \ref{thm:conv_track}, given
$\varepsilon
>0$ there exists
a constant $\bar K$ such that
%
%
\begin{equation}
\limsup_{n\tinf}\Pd \bigl\{\bigl\|\hQ_{-i_0}^n-
\cX_{-i_0}^n\bigr\| _{T^n}>\bar Kp_n \bigr\}<
\varepsilon \label{eq:SSC_stopped1}
\end{equation}
and
%
%
\begin{equation}
\limsup_{n\tinf}\Pd \bigl\{\bigl|H^n\bigr|_{T^n}+\bigl\|
\hB^n\bigr\|_{T^n}+ \bigl\|\hQ^n-\cQ^n
\bigr\|_{T^n}>\bar Kq_n \bigr\}< \varepsilon.
\label{eq:SSC_stopped2}
\end{equation}
\end{proposition}

This result is proved in the next subsection.

\begin{pf*}{Proof of Theorem \protect\ref{thm:conv_track}}
For $x\dvtx[0,u]\to\R$ denote
\[
\operatorname{Osc}\bigl(x,[s,t]\bigr):= \sup_{s\leq t_1\leq t_2\leq t}\bigl|x(t_1)-x(t_2)\bigr|.
\]
Fix $\varepsilon>0$, and let
%
%
\begin{eqnarray}\label{eq:checkOmega}
\check{\Om}^n &=&\bigl\{\bigl|H^n\bigr|_{T^n}+\bigl\|
\hB^n\bigr\|_{T^n}+\bigl\|\hQ^n-\cQ^n
\bigr\|_{T^n} \leq\bar Kq_n\bigr\} \nonumber\\
&&{}\cap\bigl\{\bigl\|
\hQ_{-i_0}^n-\cX_{-i_0}^n
\bigr\|_{T^n}\leq\bar K p_n\bigr\}
\nonumber
\\[-8pt]
\\[-8pt]
\nonumber
&&{} \cap \bigl\{\operatorname{Osc}\bigl(\theta'W^n,[s,t]
\bigr)\le p_n+L\omega(t-s), 0\le s<t\le u\bigr\}\\
&&{} \cap \bigl\{
\bigl\|W^n\bigr\|_u\leq K\bigr\}.\nonumber
\end{eqnarray}
Using Proposition \ref{thm:conv_stopped} and Lemma \ref{lem:Wtightness},
%
%
\begin{equation}
\label{84} \Pd\bigl\{\check{\Om}^n\bigr\}\geq1-4\varepsilon,
\end{equation}
provided $n$ is sufficiently large.
We begin by showing that
%
%
\begin{equation}
\label{eq:tauninfty} T^n=u\qquad \mbox{on } \check{\Om}^n.
\end{equation}

Let $M=\min_{i\in\I}(\bI\theta_{i})^{-1}$. Given $t\ge0$, we argue
that, on $\check\Om^n$, for $t<T^n$,
and all~$n$ sufficiently large,
%
%
\begin{equation}\qquad
\label{81} \theta'\hQ^n(t)\geq\eps_n:=
\bigl(1+2M^{-1}\bigr)\bar Kq_n \mbox{ implies } \sum
_{i}\tilde{B}_{ij}^n(t)=0\qquad
\mbox{for all $j$}.
\end{equation}
To see this note that if $\theta'\hQ^n(t)\geq\eps_n$, then on
$\check
{\Om}^n$,
\[
\hQ_{i_0}^n(t)\geq\cQ_{i_0}^n(t)-
\bar Kq_n=f_{i_0}^n\bigl(\theta'
\hQ ^n(t)\bigr)-\bar K q_n\geq2\bar K q_n-\bar
K q_n>0,
\]
where we use the fact that $f_{i_0}^n(x)\geq M(x\wedge\bar\kappa_n)$
for all $x\geq0$; see Remark~\ref{rem:qnproperties}. Further, on
$\check{\Om}^n$, $\theta'\hX^n(t)= \theta'\hQ^n(t)-H^n(t)\geq
\eps
_n-\bar K q_n=2M^{-1} \bar K q_n$ and $\|\hQ_{-i_0}^n-\cX_{-i_0}^n\|
_{T^n}\leq\bar K p_n$ so that, for $i\in\I_{-i_0}$,
\[
\hQ_i^n(t)\geq\check{X}_i^n(t)-
\bar K p_n=f_i^n\bigl(\theta'
\hX ^n(t)\bigr)-\bar K p_n\geq M\kappa_n-\bar
K p_n>0,
\]
where the inequalities follow from the fact that $p_n/\kappa
_n\rightarrow0$, and $f_i^n(x)\geq M\kappa_n$ for all $x\geq\kappa_n$.

Thus, we have that $\hQ_i^n(t)>0$ for all $i\in\I$ and the claim
\eqref{81} follows by Remark~\ref{rem:workconservation}.

We now argue that $G^n$ \eqref{64} remains constant when
$\sum_{i}\tilde{B}_{ij}^n(t)=0$ for all $j$. Given such $t$, let
$x\in
\X$ be defined by
$x_{ij}=-\nu_j^{-1}\tilde{B}_{ij}^n(t)$. Then $\sum_{i}x_{ij}=0$ for all
$j\in\J$. Recalling that we do not have nonbasic activities,
we conclude by the second part of Lemma \ref{lem1} that
$\frac{d}{dt}G^n(t)=-\sum_{i\sim j}\theta_i\mu_{ij}\tilde{B}_{ij}^n(t)=0$.

Combining the above argument with \eqref{81}, we have on $\check\Om^n$,
%
%
\begin{equation}
\int_{0}^{T^n}1_{(\eps_n,\infty)} \bigl(
\theta'\hQ^n(t)\bigr)\,dG^n(t)=0.
\label{eq:complementarity}
\end{equation}
By \eqref{31},
%
%
\begin{equation}
\label{82} \theta'\hQ^n(t)=\tilde W^n(t)+G^n(t),
\end{equation}
where
%
%
\begin{eqnarray}
\label{85} \tilde{W}^n(t)&=&\theta'\wh
X^n(0)+\theta'W^n(t)-R^n(t)+H^n(t),
\\
\label{96} R^n(t)&=&\sum_i
\theta_i\sum_j\eps_{ij}^n
\int_0^t\til B^n_{ij}(s)\,ds
\end{eqnarray}
and $H^n$ is as in \eqref{64}.

Thus, on $\check{\Om}^n$ and for $t\in[0,T^n]$, the
triplet $(\theta'\hQ^n,\tilde{W}^n,G^n)$ satisfies the following relations:
namely, $G^n(0)=0$; as shown in the proof of Proposition \ref{prop2},
$G^n$ is nondecreasing; and since $\hQ^n$ takes values in the positive
orthant, $\theta'\hQ^n\geq0$.
With these properties, along with \eqref{eq:complementarity} and
\eqref
{82}, we are
in a position to apply the oscillation inequalities from \cite{wilSRBM}, Theorem
5.1,
to conclude that
%
%
\begin{equation}
\operatorname{Osc}\bigl(\theta'\hQ^n,[s,t]\bigr)+ \operatorname{Osc}
\bigl(G^n,[s,t]\bigr)\leq4 \bigl(\operatorname{Osc}\bigl(\tilde{W}^n,[s,t]
\bigr)+\eps_n\bigr) \label{eq:oscillation1}
\end{equation}
for all $0\leq s\leq t\leq T^n$.
The precise constant in the inequality follows the proof of the
one-dimensional case \cite{wilSRBM}, page 15.

To prove \eqref{eq:tauninfty} we use the oscillation inequality to show
first that $T^n=\tau^n$. We then show that $\tau^n= u$.
By the definition of $T^n$, using \eqref{32}, we have
%
%
\begin{equation}
\label{83} \bigl\|R^n\bigr\|_{T^n}\leq J\|\theta\|
\varepsilon_M^n\varrho^n\leq p_n,
\end{equation}
and using the definition of $\check\Om^n$, we have on $\check\Om^n$,
\[
\bigl|H^n\bigr|_{T^n}\leq\bar Kq_n.
\]
By Lemma \ref{lem:Wtightness}, $\operatorname{Osc}(\theta'W^n,[s,t])\leq
p_n+L\omega(t-s)$.
Hence on $\check\Om^n$,
%
%
\begin{equation}
\operatorname{Osc}\bigl(\tilde{W}^n,[s,t]\bigr)\leq (2+\bar
K)q_n+L\omega(t-s). \label{eq:tildeWosc}
\end{equation}
Using \eqref{eq:oscillation1}, we then have
\[
\operatorname{Osc}\bigl(G^n,[s,t]\bigr)\leq4(2+\bar
K)q_n+4L\omega(t-s))+4\eps_n
\]
for $s<t\leq T^n$.
Recalling that $r_n/q_n\tinf$, we then must have $\zeta^n\geq\tau^n$
on $\check{\Om}^n$. Further, since $\|\hB^n\|_{T^n}\leq\bar{K}q_n$ on
$\check{\Om}^n$
and, recalling that $q_n/s_n\rightarrow0$, clearly $\sigma^n\geq
T^n$. We
conclude that $T^n=\tau^n$ on $\check{\Om}^n$.

To prove \eqref{eq:tauninfty} it then remains to show that $\tau^n= u$.
Following the same arguments leading to \eqref{68}, on the event
that $\{\tau^n<u\}$ one has $\theta'\wh Q^n(\tau^n)\ge c(\varrho_n-K)$,
for some constant
$c>0$. Recall that $\varrho_n\to\iy$. On the other hand,
using \eqref{82}, \eqref{eq:oscillation1}, and the bounds on $R^n$
and $H^n$,
%
%
\begin{equation}
\bigl|\theta' \hQ^n\bigr|_{\tau^n}\leq\bigl|
\tilde{W}^n\bigr|_{\tau^n}+\bigl|G^n\bigr|_{\tau^n} \le5
\bigl(\|\theta\|K+(1+\bar K)q_n\bigr)+4\eps_n
\label{eq:interim3}
\end{equation}
on $\check{\Om}^n(\varepsilon)$. Thus $|\theta' \hQ^n|_{\tau^n}$ is
bounded on this
event. This shows that $\tau^n=u$ on $\check\Om^n$
for all sufficiently large $n$.
We have thus proved that $T^n=\tau^n=u$ on $\check{\Om}^n$, for
large~$n$.

Since $\eps$ is arbitrary, using Proposition \ref{thm:conv_stopped},
\eqref{84} and \eqref{83}, we obtain
%
%
\begin{equation}
\bigl|H^n\bigr|_{u}+\bigl\|\hB^n\bigr\|_{u}+\bigl\|
\hQ^n-\cQ^n\bigr\|_{u}+\bigl|R^n\bigr|_{u}
\to0\qquad\mbox{in probability.} \label{eq:convu}
\end{equation}

Using Lemma \ref{lem2}(iii) and arguing along the lines of step 3 of
the proof of
Proposition \ref{prop2} (showing, in particular, that $\bar T^n\to
\bar T$ in probability)
gives $W^n\To W$. Using \eqref{85}, \eqref{eq:convu} and recalling that
we assumed zero initial conditions, we have $\tilde W^n\To\theta'W$.

We conclude that the process $\theta'\wh Q^n$ satisfies
$\theta'\wh Q^n=\tilde W^n+G^n\ge0$, where $\tilde W^n\To\theta'W$,
and $G^n$ is nondecreasing and satisfies, with probability arbitrarily
close to one,
\[
\int_0^\iy1_{(\eps_n,\iy)}\bigl(
\theta'\wh Q^n(t)\bigr)\,dG^n(t)=0.
\]
It is a standard fact that these properties suffice
to characterize the limit behavior, namely
that $\theta'\wh Q^n\To\Gam[\theta'W]=Q^*$; see, for example, the proof
of \cite{wilSRBM}, Theorem 4.1.
By \eqref{86}, the uniform convergence of $f^n$ to $f$, \eqref
{eq:qconv} and the continuous mapping theorem, $\check Q^n\To q(Q^*)$.
Hence by \eqref{eq:convu}, $\wh Q^n\To q(Q^*)$, thus
$C(\wh Q^n)\To C(q(Q^*))=C_*(Q^*)$.
Another application of the continuous mapping theorem gives
\[
\int_{0}^u C\bigl(\hQ^n(t)\bigr)\,dt
\Rightarrow\int_{0}^u C_{*}
\bigl(Q^*(t)\bigr)\,dt,
\]
which completes the proof the theorem.
\end{pf*}

\subsection{\texorpdfstring{Proof of Proposition \protect\ref{thm:conv_stopped}}
{Proof of Proposition 5.1}}

The key idea in the proof is to identify an event occurring with high
probability on which the policy\vspace*{1pt} self tunes the balance between $\cX^n$
and $\hQ^n$: when the process $\hX^n$ goes ``out of balance,'' namely,
when $\|\hQ_{-i_0}^n(t)-\cX_{-i_0}^n(t)\|>cp_n$, the occupancy process
$\hB^n$ re-adjusts quickly so as to pull the process $\hQ^n$ back
toward $\cX^n$.

Throughout the remainder of the analysis we fix $\varepsilon>0$; $T^n$ is
as in \eqref{eq:checktaudefin}. Define
%
%
\begin{eqnarray}\label{eq:omegatilde2}
\qquad\quad \Om_1^n&=&\bigl\{ \bigl|\theta'
\hX^n\bigr|_{T^n}\leq K\bigr\},
\nonumber
\\[-8pt]
\\[-8pt]
\nonumber
\Om_2^n&=&\Bigl\{ \max
_{\La^n=\wh A^n,V^n,W^n}\bigl\|\La^n(t)-\La^n(s)\bigr\|\leq
p_n+L\omega(t-s), 0\leq s\leq t\leq u\Bigr\},
\end{eqnarray}
where, with an abuse of notation, $K=K(\eps)$ and $L=L(\eps)$ will be
chosen (possibly)
larger than the values from Lemma \ref{lem:Wtightness}. Let
%
%
\begin{equation}
\label{87} \Om^n=\Om_1^n\cap
\Om_2^n.
\end{equation}

\begin{lemma}\label{lem:omegatilde} Suppose that the assumptions of Theorem \ref
{thm:conv_track} hold. Then
$K$ and $L$ can be chosen so that $\Pd\{\Om^n\}\geq1-\varepsilon$.
Moreover, on $\Om^n$, and for $0\leq s\leq t\leq T^n$,
%
%
\begin{equation}
\label{eq:omegatilde3} \bigl\|\cX^n(t)-\cX^n(s)\bigr\|\leq c\bigl(
\min\bigl[n^{1/2}s_n(t-s)+p_n,r_n
\bigr]+\omega(t-s)\bigr),
\end{equation}
where $c$ is a constant not depending on $n,s,t$.
\end{lemma}
Lemma \ref{lem:omegatilde} is proved in Section~\ref{sec:lem}.

Throughout what follows, $K$ and $L$ are as in the above lemma, and fixed,
$\Om^n$ is as in \eqref{87} and $T^n$ as in \eqref{eq:checktaudefin}.
Given strictly positive constants $\{c_k^1, k\in\I\}$ define the
following times:
%
%
\begin{eqnarray}
\tau_{2,k}^n&:=&\inf\bigl\{s\geq0\dvtx\bigl|
\Del_k^n(s)\bigr| >c_{k}^1p_n
\bigr\}\wedge T^n,\label{eq:tau2kdefin}
\\
\tau_{1,k}^n&:=&\sup\bigl\{s\leq\tau_{2,k}^n
\dvtx\bigl|\Del_k^n(s)\bigr|\leq c_{k}^1p_n/2
\bigr\}, \label{eq:tau1kdefin}
\end{eqnarray}
where, throughout,
\[
\Del^n_k=\wh Q^n_k-\check
X^n_k.
\]
Let
\[
\Om^n_{k,U}=\Om^n\cap\bigl\{
\tau^n_{2,k}<T^n\bigr\}\cap \bigl\{
\Del_k^n\bigl(\tau^n_{2,k}\bigr) >
c_{k}^1p_n\bigr\},
\]
where $U$ is mnemonic for ``up,'' and define analogously $\Om^n_{k,D}$,
with ``$> c_{k}^1p_n$'' replaced by ``$< -c_{k}^1p_n$,''
where $D$ is mnemonic for ``down.''
Note that the jumps of $\wh Q^n$ and $\check X^n$ are of order
$n^{-1/2}$ while
$p_n\gg n^{-1/2}$. Moreover, the initial condition is assumed to be zero.
As a result, we have on the event $\Om^n_{k,U}$ (resp., $\Om^n_{k,D}$)
that $\tau^n_{1,k}\in[0,\tau^n_{2,k} )$,
and
%
%
\begin{equation}
\label{91} \Del^n_k(s)\ge c^1_kp_n/2\qquad
\bigl(\mbox{resp., } \le-c^1_kp_n/2\bigr)
\mbox{ for all } \tau^n_{1,k}\le s<\tau^n_{2,k}.
\end{equation}

The proof of Proposition \ref{thm:conv_stopped} will be based on
showing that
$\tau_{2,k}^n\geq T^n$ on $\Om^n$ for all $k\in\I_{-i_0}$. This
statement is proved inductively.

\begin{proposition}\label{prop:localadjustment1} Suppose that the assumptions of Proposition \ref
{thm:conv_stopped} hold.
Then the following holds on the event $\Om^n$.
Let $k\in\I_{-i_0}$ be either $1$ or such that, for all $l<k$,
%
%
\begin{equation}
\bigl|\Del_l^n\bigr|_{T^n}\leq c_l^1p_n
\label{eq:inductionarg}
\end{equation}
for some constants $c_l^1$, $l<k$.
Then there exists a constant $c_k^1$, such that if $\tau^n_{2,k}$
\eqref{eq:tau2kdefin} is defined with $c_k^1$, then $\tau_{2,k}^n\geq
T^n$. Consequently, there exist constants $c_k^1, k\in\I_{-i_0}$
such that, for all $k\in\I_{-i_0}$,
%
%
\begin{equation}
\bigl|\Del_k^n\bigr|_{T^n}\leq c_k^1p_n.
\label{eq:inductionarg1}
\end{equation}
\end{proposition}

Recall \eqref{eq:gndefin} and assume, without loss of generality, that
$q_n=n^{\delta}p_n$ for some $\delta>0$ such that $q_n\ll r_n$.
Recall that $\wh B^n_{ij}=B^n_{ij}=0$ if $i\nsim j$.

The next proposition is where the perturbation $f^n$ of the function
$f$ is used in an important way.

\begin{proposition}\label{prop:localadjustment2} Suppose that the assumptions of Proposition \ref
{thm:conv_stopped} hold. Then, there exists a constant $\gamma$ such
that, on the event $\Om^n$,
\[
\bigl\|\hB^n\bigr\|_{T^n}\leq\gamma q_n.
\]
\end{proposition}

The proofs of Propositions \ref{prop:localadjustment1} and \ref
{prop:localadjustment2} appear in Sections~\ref{sec:prop1proof} and
\ref{sec:prop2proof}, respectively. Henceforth we let $M_x$ be a
uniform Lipschitz constant of $f^n$ on $[0,x]$. Throughout the proofs
we use $c_1,c_2,\ldots$ to denote strictly positive constants that do not
depend on $n$.

\begin{pf*}{Proof of Proposition \protect\ref{thm:conv_stopped}}
Equation \eqref{eq:SSC_stopped1} follows directly from Proposition
\ref
{prop:localadjustment1} and the definition of $\Omega^n$.

To prove \eqref{eq:SSC_stopped2}, we will first prove that
%
%
\begin{equation}
\bigl\|\hQ^n-\cX^n\bigr\|_{T^n}\leq c_1
q_n\label{eq:22}.
\end{equation}

To this end, by Proposition \ref{prop:localadjustment2} and identity
\eqref{06} we have that
%
%
\begin{equation}
\bigl\|\hX^n-\hQ^n\bigr\|_{T^n}\leq\bigl\|\hB^n\bigr\|
_{T^n}\leq c_2 q_n.\label{eq:21}
\end{equation}
By Proposition \ref
{prop:localadjustment1} and the fact that $p_n\ll q_n$ we have that
%
%
\begin{equation}
\bigl\|\hQ_{-i_0}^n-\cX_{-i_0}^n
\bigr\|_{T^n}\leq c_3q_n.\label{eq:20}
\end{equation}

By definition $\theta'\hX^n=\theta'\cX^n$ whenever $\theta'\hX
^n\geq0$
and $\theta
\cX^n=0$ otherwise. Using~\eqref{eq:21} and the nonnegativity of $\hQ
^n$, we get that $\theta'\hX^n(t)\geq-c_4 q_n$ for all $t\leq T^n$
and, in turn, that $|\theta'\hX^n-\theta'\cX^n|_{T^n}\leq c_5 q_n$. Thus
\[
\bigl|\hX_{i_0}^n-\cX_{i_0}^n\bigr|_{T^n}
\leq\bigl|\theta'\hX^n-\theta'\cX
^n\bigr|_{T^n}+\bigl\|\hX_{-i_0}^n-
\cX_{-i_0}^n\bigr\|_{T^n}\leq c_6
q_n.
\]
Applying \eqref{eq:21} we conclude that
$|\hQ_{i_0}^n-\cX_{i_0}^n|_{T^n}\leq c_7 q_n$ and, together with
\eqref
{eq:20}, that \eqref{eq:22} holds.

Next, since $|\theta'\hX^n|_{T^n}\leq K$ on $\Om^n$, we have
$|\theta
'\hQ^n|_{T^n}\leq
|\theta'\hX^n|_{T^n}+|H^n|_{T^n}\leq2K$.
Hence, by Assumption \ref{asum:continuousselection} and the definition
of $\check Q^n$ and $\check X^n$,
\[
\bigl\|\cX^n-\cQ^n\bigr\|_{T^n}\leq M_{2K}\bigl|
\theta'\hQ^n-\theta'\hX
^n\bigr|_{T^n}\leq c_8 q_n,
\]
where we used \eqref{eq:21}.
Using this along with \eqref{eq:22}, we have on $\Om^n$,
%
%
\begin{equation}
\bigl\|\hQ^n-\cQ^n\bigr\|_{T^n}\leq c_{9}
q_n. \label{eq:interim55}
\end{equation}

Finally, recall that $-H^n=\sum_{i,j}\theta_i\hB_{ij}^n$ so that
%
%
\begin{equation}
\label{90} \bigl|H^n\bigr|_{T^n}\le c_{10} q_n,
\end{equation}
by Proposition \ref{prop:localadjustment2}. Combining \eqref
{eq:interim55}, \eqref{90} and Proposition \ref{prop:localadjustment2}
we conclude that for any $\eps$ there exists
a constant $c_{11}$ independent of $n$, such that
\[
\Pd\bigl\{\bigl|H^n\bigr|_{T^n}+\bigl\|\hB^n\bigr\|_{T^n}+
\bigl\|\hQ^n-\cQ^n\bigr\|_{T^n}> c_{11}q_n
\bigr\} <\varepsilon.
\]
This proves \eqref{eq:SSC_stopped2}.
\end{pf*}

\subsection{\texorpdfstring{Proof of Proposition \protect\ref{prop:localadjustment1}}
{Proof of Proposition 5.2}}\label{sec:prop1proof}
We begin by stating a sequence of lemmas that provide estimates on various
properties of the dynamics.
They are proved in Section~\ref{sec:lem}, along with Lemma \ref
{lem:omegatilde} above.
Throughout this and the next subsection, $\varepsilon$ is fixed, and the
assumptions
of Theorem \ref{thm:conv_track} are in force. Moreover, the statements
of the lemmas
are understood to be \textit{on the event $\Om^n$.}
For $f\dvtx[0,\iy)\to\R^k$ and $0\le s\le t$, denote $f[s,t]:=f(t)-f(s)$.

The proposition provides estimates concerning the r.v.'s $\tau^n_{2,k}$,
that are based on the lemmas below. At the same time, the proof of the
proposition involves choosing the constants $c^1_k$,
used to define these r.v.'s.
It will therefore be important to specify which of the estimates,
stated in the lemmas
(at least those that involve~$\tau^n_{2,k}$), depend on $c^1_k$, and
which do not.

Define the processes
%
%
\begin{equation}
\label{92} \mB_{ij}^n(t):=\hB_{ij}^n(t)+
\mu_{ij}^n\int_0^t \hB
_{ij}^n(s)\,ds,\qquad i\in\calI,j\in\calJ.
\end{equation}

\begin{lemma}\label{lem:aux3}
Fix $k\in\I_{-i_0}$. Suppose that $|\Del_l^n|_{T^n}\leq c_l^1p_n$ for
all $l<k$ for some constants $c_l^1$, $l<k$. Then there exists a
constant $\gamma_1$, not depending on
the constant $c^1_k$ with which $\tau^n_{2,k}$ is defined,
such that, for all $0\leq s\leq t\leq T^n$,
%
%
\begin{equation}
\bigl|\mB_{kj}^n[s,t]\bigr|\leq \gamma_1\bigl(
n^{1/2}s_n(t-s)+p_n+\omega(t-s)\bigr),\qquad j\in
\calJ(k). \label{eq:mBincrements}
\end{equation}
Moreover, there exists a constant $\gamma_2$ (that may depend on $c^1_k$)
such that \eqref{eq:mBincrements}
holds for $j=\jbar(k)$ and $0\leq s\leq t\leq\tau_{2,k}^n$, with
$\gamma_1$ replaced by
$\gamma_2$.
\end{lemma}
%
\begin{lemma}\label{lem:qtightnessuptostop}
Fix $k\in\I_{-i_0}$. Then there exists a constant $\gamma$ (that may
depend on
$c^1_k$) such that for all $0\leq s\leq t\leq\tau_{2,k}^n$,
\[
\bigl|\hQ_k^n[s,t]\bigr|\leq\gamma\bigl( n^{1/2}s_n(t-s)+p_n+
\omega(t-s)\bigr).
\]
\end{lemma}

To state the last preliminary lemma, let
%
%
\begin{equation}
\mathcal{L}(k)=\bigl\{i\in\I\dvtx i\leq k, i\sim\jbar(k)\bigr\}
\label{eq:lkdefin}
\end{equation}
be the set of customer classes that are not higher than $k$ in the
hierarchy and
are connected to the parent node $\jbar(k)$. Recall that if $\tau
^n_{2,k}<T_n$,
then one of the two events $\Om^n_{k,U}$, $\Om^n_{k,D}$ must occur,
and consequently, one of the two inequalities specified in \eqref{91} holds.
In the former case, any
service completion in pool $\jbar(k)$ during $[\tau_{1,k}^n,\tau
_{2,k}^n)$ is followed,
under our tracking policy, with an admission of a customer from one of
the queues in the set $\mL(k)$. In the latter case, no class-$k$
customers are admitted to pool $\jbar(k)$ on
$[\tau_{1,k}^n,\tau_{2,k}^n)$. The following lemma is based on these
two properties.

\begin{lemma}\label{lem:aux2}
Fix $k\in\I_{-i_0}$, and let $j=\jbar(k)$.
Then there exists a constant $\gamma>0$, not depending on $c^1_k$,
such that the following holds:

\begin{longlist}[(ii)]
\item[(i)]
On $\Om^n_{k,U}$ one has
%
%
\begin{eqnarray}\label{eq:aux2_1}
\sum_{l\in\mathcal{L}(k)}\mB_{lj}^n[s,t]
\geq\gamma n^{1/2}(t-s)-L\omega(t-s)-p_n,
\nonumber
\\[-8pt]
\\[-8pt]
\eqntext{\tau_{1,k}^n\leq s\leq t\leq\tau_{2,k}^n.}
\end{eqnarray}
\item[(ii)]
On $\Om^n_{k,D}$ one has
%
%
\begin{equation}
\qquad\mB_{kj}^n[s,t] \leq-\gamma n^{1/2}(t-s)+L
\omega(t-s)+p_n,\qquad \tau_{1,k}^n\leq s\leq t\leq
\tau_{2,k}^n. \label{eq:aux2_2}
\end{equation}
\end{longlist}
\end{lemma}

\begin{pf*}{Proof of Proposition \protect\ref{prop:localadjustment1}}
Note that all statements regard the event $\Om^n$~\eqref{87}.

It is required to show that there exists a constant
$c^1_k$, with which $\tau^n_{2,k}$ is defined, such that $\tau
^n_{2,k}\ge T_n$.

To this end, let us analyze the event $\Om^n\cap\{\tau_{2,k}^n<T^n\}$,
considering
separately the two sub-events $\Om^n_{k,U}$ and $\Om^n_{k,D}$.
The goal is to show that one can choose $c^1_k$ so that
the two events are empty, provided $n$ is sufficiently large.

We start with the former. Fix $k\in\calI$ [and note that in the case $k=1$,
the set $\mL(k)$ is simply $\{1\}$],
and denote $j=\jbar(k)$.
The goal of showing that $\Om^n_{k,U}$ is empty (for suitable $c^1_k$
and large $n$) is achieved
by arguing that there exists a constant $c$, not depending on $c^1_k$
or $n$,
such that on $\Om^n_{k,U}$,
%
%
\begin{equation}
\label{99} \Del^n_k(t)\le\Del^n_k
\bigl(\tau^n_{1,k}\bigr)+cp_n \qquad\mbox{for all } t
\in \bigl[\tau^n_{1,k},\tau^n_{2,k}\bigr).
\end{equation}

To this end, consider any $\tau^n_{1,k}\le s\le t < \tau^n_{2,k}$.
Using \eqref{31} and \eqref{92} we have
%
%
\begin{eqnarray}\label{eq:interim77}
\sum_{l\in\mathcal{L}(k)}\hQ_l^n[s,t]&=&
\nonumber
\sum_{l\in\mathcal{L}(k)}W_l^n[s,t]
\\
&&{}-\sum_{l\in\mathcal{L}(k)}\sum
_{m\in\calJ(l)}\mB_{lm}^n[s,t]
\\
&&{}- \sum_{l\in\mathcal{L}(k)} \mB_{lj}^n[s,t],
\nonumber
\end{eqnarray}
where we used the fact that $\jbar(l)=j$ for any $l\in\calL(k)$.
By the assumption of the proposition, $|\Del_l^n|_{T^n}\le c^1_lp_n$, $l<k$,
by which $\tau^n_l=T^n$ for $l<k$. Thus in view of Lemma \ref{lem:aux3},
estimate \eqref{eq:mBincrements} is valid for all $l<k$, and for all
$0\le s\le t\le T^n$.
Moreover, the constants in this estimate do depend on $c^1_l$, $l<k$,
but not on $c^1_k$.
As a result,
%
%
\begin{equation}
\sum_{l\in\mathcal{L}(k)}\sum_{m\in\calJ(l)} \bigl|
\mB_{lm}^n[s,t]\bigr| \leq c_1\bigl(
n^{1/2}s_n(t-s)+p_n+\omega(t-s)\bigr),
\label{eq:interim88}
\end{equation}
where $c_1$ does not depend on $c^1_k$.
We bound the first line of \eqref{eq:interim77} using \eqref
{eq:omegatilde2}, the second line using \eqref{eq:interim88} and the
third line using Lemma \ref{lem:aux2}. Here we recall that the
constants in Lemma \ref{lem:aux2} do not depend on $c_k^1$. In turn,
we have
%
%
\begin{eqnarray}
\sum_{l\in\mathcal{L}(k)}\hQ_l^n[s,t]&
\leq& c_2\bigl(p_n+\omega(t-s)+ n^{1/2}s_n(t-s)
\bigr)-c_3 n^{1/2}(t-s)
\nonumber
\\[-8pt]
\\[-8pt]
\nonumber
& \leq& -c_4 n^{1/2}(t-s)+c_5
\bigl(p_n+\omega(t-s)\bigr), \label{eq:1}
\end{eqnarray}
or positive constants $c_2, \ldots,c_5$ that do not depend on $c^1_k$.
In the second inequality above we used the fact that $s_n\rightarrow0$.
Thus, using Lemma \ref{lem:qtightnessuptostop}
\[
\sum_{l<k}\bigl|\hQ_l^n[s,t]\bigr|
\leq c_8 \bigl( n^{1/2}s_n(t-s)+p_n+
\omega (t-s) \bigr).
\]
Above, the constant may depend on $c^1_l$, $l<k$, but not on $c^1_k$.
Thus by \eqref{eq:1} and using the fact that $s_n\rightarrow0$ and
applying Lemma \ref{lem:omegatilde} to bound $|\cX_k^n[s,t]|$, we obtain
\begin{eqnarray*}
\Del^n_k(t)=\hQ_k^n(t)-
\cX_k^n(t)&\leq&\Del^n_k(s)
-c_{11} n^{1/2}(t-s)+c_{12}\bigl(p_n+
\omega(t-s)\bigr).
\end{eqnarray*}
Note that we have
%
%
\begin{equation}
-c_{11} n^{1/2}(t-s)+c_{12}\omega(t-s)\leq
c_{13} p_n\label {eq:interim99}
\end{equation}
for all sufficiently large $n$. Indeed, the function $\nu\to-c_{11}
n^{1/2}\nu+c_{12}\omega(\nu)$ is concave and it is easily verified that
the unique maximum is bounded by $c_{13} n^{-{\alpha_{\omega
}}/{(2(1-\alpha_{\omega}))}}$. Since $\alpha_{\omega}\in(1/3,1/2)$
the maximum
is further bounded by $c_{14} n^{-1/4}\ll p_n$. Thus, we conclude that
\eqref{eq:interim99} holds.
Choosing $c_k^1$ sufficiently large, we then have that $\tau
_{2,k}^n\geq T^n$ must hold on the event
$\Om^n_{k,U}$. In other words, $\Om^n_{k,U}$ is empty.

Next we consider the event $\Om^n_{k,D}$.
Arguing as above, using the second part of Lemma \ref{lem:aux2}, we have
\[
\hQ_k^n(t)\geq\hQ_k^n(s)-c_{15}
\bigl(\omega(t-s)+ n^{1/2}s_n(t-s)\bigr) +c_{16}
n^{1/2}(t-s).
\]
Bounding $\cX_k^n[s,t]$ and using \eqref{eq:omegatilde3} as before we have
\[
\Del_k^n(t)\geq \Del_k^n(s)-c_{17}
\bigl(\omega(t-s)+p_n\bigr)+c_{18} n^{1/2}(t-s).
\]
Similarly to the above, $\Del_k^n(t)\leq\Del_k^n(\tau_{1,k}^n)+c_{19}p_n$
for all $t\in[\tau_{1,k}^n,\tau_{2,k}^n)$.
Re-choosing $c_{k}^1\geq4c_{20}$ we then conclude that $\tau
_{2,k}^n\geq T^n$
on $\Om^n_{k,D}$.

This proves the first part of the proposition. The second part is
argued inductively using the above. If \eqref{eq:inductionarg} holds
for all $l<k$ (or for the induction basis $k=1$) and since we proved
that $\tau^n_{2,k}\geq T^n$, the definition of $\tau_{2,k}^n$ implies
that $|\Del^n_k(s)|\le c_k^1p_n$ for all $s<T_n$. It is not hard to see
that the jumps of both $\wh Q^n$ and $\check X^n$ are $O(n^{-1/2})$.
Recalling that $p_n \gg n^{-1/2}$, it follows that \eqref
{eq:inductionarg} holds for $k$ (with a suitable constant $c_k^1$). We
conclude that there exist constants, that with abuse of notation we
still denote by $\{c_k^1, k\in\I\}$, such that, on $\Om^n$, for all
$k\in\I_{-i_0}$,
$|\Del_k^n|_{T^n}\leq{c}_k^1p_n$. This concludes the proof.
\end{pf*}

\subsection{\texorpdfstring{Proof of Proposition \protect\ref{prop:localadjustment2}}
{Proof of Proposition 5.3}}\label
{sec:prop2proof}
We begin by stating a sequence of auxiliary lemmas that are proved in
Section~\ref{sec:lem}. As before, the statements of the lemmas are
understood to be \textit{on the event $\Om^n$} and the assumptions
of Theorem \ref{thm:conv_track} are in force. Fixing throughout
$\delta
$ such that $q_n=p_nn^{\delta}\ll r_n$ we let $\vartheta
_n=n^{-(1/2-\delta)}$. Below the constants $z_j^*,j\in\J$, are as in
\eqref{eq:zdefin}.

The following relates the process $\mB^n$ to the idleness process.

\begin{lemma}\label{lem:idleness}
Fix $j\in\J$ with $j\ne j_0$. Suppose that there exists a constant
$\gamma_1$ such that
\[
\sum_{k\in\calI(j)}\bigl|\mB_{kj}^n[s,t]\bigr|
\leq \gamma_1\bigl( n^{1/2}s_n(t-s)+p_n+
\omega(t-s)\bigr)
\]
for all $0\leq s\leq t\leq T^n$. Then there exists a constant $\gamma
_2$ such that $|\hI_j^n|_{T^n}\leq\gamma_2 p_n$.

Consequently, if $|\Del_k^n|_{T^n}\leq c_k^1p_n$ for $k\in\I_{-i_0}$
and constants $c_k^1$, $k\in\I_{-i_0}$, then there exists a constant
$\gamma_2$ such that $|\hI_j^n|_{T^n}\leq\gamma_2 p_n$, for all
$j\neq j_0$.
\end{lemma}

\begin{lemma}\label{lem:XtoI} There exists a constant $\gamma$ so that for all
$s,t\leq
T^n$ with $|t-s|\leq\vartheta_n$,
\[
\biggl|\theta'\hX^n[s,t]-z_{j_0}^*
n^{1/2}\int_s^t \hI
_{j_0}^n(u)\,du \biggr|\leq\gamma q_n
\]
and
\[
\theta'\hX^n[s,t]\geq-\gamma p_n.
\]
\end{lemma}

\begin{lemma}\label{lem:cXoscillation} Fix $k\in\I_{-i_0}$. Then there
exists a constant $\gamma$ such that for all $s,t\leq T^n$ with
$|t-s|\leq\vartheta_n$
\[
\bigl|\cX_k^n[s,t]\bigr|\leq\gamma q_n.
\]
\end{lemma}

\begin{lemma}\label{lem:Bbound} Fix $k\in\I_{-i_0}$ and $j\neq j_0$.
Then there exists a constant $\gamma$ such that
\[
\bigl|\hB_{kj}^n\bigr|_{T^n} \leq\gamma q_n.
\]
\end{lemma}

\begin{lemma} \label{lem:oscillationI} There exists a constant $\gamma$
such that the following holds. For $s,t\leq T^n$ with $|t-s|\leq
\vartheta_n$ and such that $\hI_{j_0}^n(u)>0$ for all $u\in[s,t)$,
we have
\[
\hI_{j_0}^n[s,t]\leq-\cX_{i_0}^n[s,t]+
\gamma q_n\leq -\frac{1}{2} z_{j_0}^*\int
_{s}^tI_{j_0}^n(u)\,du+
\gamma q_n.
\]
\end{lemma}

\begin{pf*}{Proof of Proposition \protect\ref{prop:localadjustment2}}
By Lemma \ref{lem:Bbound}, $\sum_{j\neq j_0}\sum_{i\neq i_0}|\hB
_{ij}^n|_{T^n}\leq c_1 q_n$. By Proposition \ref{prop:localadjustment1}
and Lemma \ref{lem:idleness} we have that $\sum_{j\neq j_0}|\hI
_j^n|_{T^n}\leq c_2 p_n$. Thus, using the identity $\hI_{j}^n=-\sum_{i}\hB_{ij}^n$ and since $p_n\ll q_n$, we have
\[
\sum_{j\neq j_0}\bigl|\hB_{i_0j}^n\bigr|_{T^n}
\leq\sum_{j\neq j_0}\bigl|\hI _j^n\bigr|_{T^n}+
\sum_{j\neq j_0}\sum_{i\in\I(j)}\bigl|
\hB_{ij}^n\bigr|\leq c_3 q_n.
\]

To prove the proposition it only remains to show that
%
%
\begin{equation}
\bigl|\hI _{j_0}^n\bigr|_{T^n}\leq c_4
q_n\label{eq:j0bound},
\end{equation}
in which case we will
have by the same argument that $|\hB_{i_0j_0}^n|_{T^n}\leq c_5 q_n$.
Together with Lemma \ref{lem:Bbound}, this would allow us to conclude
that $\|\hB^n\|_{T^n}\leq c_6 q_n$ as required.

The remainder of the argument is dedicated to the proof of \eqref
{eq:j0bound}. To that end, fix $\zeta>2\gamma$ with $\gamma$ as in
Lemma \ref{lem:oscillationI}, and let
\[
\tau_1^n=\inf\bigl\{t\geq0\dvtx\hI_{j_0}^n(t)>
2\zeta q_n\bigr\}\wedge T^n
\]
and
\[
\tau_0^n=\sup\bigl\{t\leq\tau_1^n
\dvtx\hI_{j_0}^n(t)\leq\zeta q_n\bigr\}\wedge
T^n.
\]
Argue by contradiction and assume that $\tau_1^n<T^n$ and, in
particular, $\tau_0^n<T^n$. Consider the interval $[\tau_0^n,(\tau
_0^n+\vartheta_n)\wedge\tau_1^n)$. By Lemma \ref{lem:oscillationI} it
holds, for $s,t\in[\tau_0^n,(\tau_0^n+\vartheta_n)\wedge\tau
_1^n)$ that
%
%
\begin{equation}
\hI_{j_0}^n[s,t]\leq -\frac{1}{2}
z_{j_0}^*n^{1/2}\int_{s}^t
\hI_{j_0}^n(u)\,du+\gamma q_n. \label{eq:Idec}
\end{equation}
In particular, $\hI_j^n(u)\leq\zeta q_n+\gamma q_n$ for all $u\in
[\tau
_0^n,(\tau_0^n+\vartheta_n)\wedge\tau_1^n)$, and it must be the case
that $\tau_1^n\geq\tau_0^n+\vartheta_n$. Since $\hI_{j_0}^n\geq
\zeta
q_n$ on $[\tau_0^n,\tau_1^n)$ and $n^{-1/2}\ll\vartheta_n$ we also
have by \eqref{eq:Idec} that
\[
\hI_{j_0}^n\bigl(\tau_0^n+4
\bigl(z_{j_0}^*\bigr)^{-1}n^{-1/2}\bigr)\leq2\zeta
q_n-2\zeta q_n +\gamma q_n\leq\gamma
q_n.
\]
Since $\tau_0^n+4(z_{j_0}^*)^{-1}n^{-1/2}\in(\tau_0^n,\tau_1^n)$ this
contradicts the definition of $\tau_0^n$. We conclude that
$\tau_1^n\geq T^n$ and, since the jumps of $\hI^n$ are of size
$O(n^{-1/2})$ and $n^{-1/2}\ll q_n$, that $|\hI_{j_0}^n|_{T^n}\leq
3\zeta q_n$. This establishes \eqref{eq:j0bound}\vspace*{1pt} and completes the
proof of the proposition.
\end{pf*}

\subsection{Proofs of auxiliary lemmas}\label{sec:lem}
\mbox{}
\begin{pf*}{Proof of Lemma \protect\ref{lem:omegatilde}}
It follows directly from Lemma \ref{lem:Wtightness} that
$\Pd\{\Om_2^n\}\geq1-\varepsilon/2$ for a sufficiently large $L$.
To treat $\Om_1^n$, recall that
%
%
\begin{equation}
\label{94} \theta'\hX^n(t)=\theta'
\hX^n(0)+\theta'W^n(t)+G^n(t)-R^n(t).
\end{equation}
By the definition of $T^n$, it is easy to see that $|R^n|_{T^n}$ and
$|G^n|_{T^n}$
are uniformly bounded.
Hence using Lemma \ref{lem:Wtightness}, $K$ can be chosen so that
$\Pd\{\Om_1^n\}\geq1-\varepsilon/2$. This shows the first assertion of
the lemma.

By \eqref{19}
\[
\bigl\|\hX^n[s,t]\bigr\|\leq \bigl\|W^n[s,t]\bigr\|+\sum
_{i,j} n^{1/2}\bigl(\mu_{ij}+
\varepsilon_{ij}^n\bigr)\int_s^t
\bigl|\hB _{ij}^n(s)\bigr|\,ds.
\]
By the definition of $\Om^n_2$ and the fact that
$\|\hat{B}^n\|_{T^n}\leq s_n$ and using Assumption \ref
{asum:continuousselection} to write
%
%
\begin{equation}
\label{95} \bigl\|\cX^n[s,t]\bigr\|\leq M_K\bigl|\theta'
\hX^n[s,t]\bigr|,
\end{equation}
it follows that
\[
\bigl\|\check X^n[s,t]\bigr\|\le c\bigl(\sqrt n s_n(t-s)+p_n+L
\omega(t-s)\bigr).
\]
To prove the result, it remains to show that
%
%
\begin{equation}
\label{93}\bigl \|\check X^n[s,t]\bigr\|\le c\bigl(r_n+L\omega(t-s)
\bigr).
\end{equation}
We use \eqref{94} and \eqref{95}. The increment of $W^n$ can be bounded
as before,
while that of $G^n$ is bounded using the definition of $T^n$,
specifically $\zeta^n$.
Thus
\[
\bigl\|\check X^n[s,t]\bigr\|\le c\bigl(r_n+L\omega(t-s)
\bigr)+\bigl|R^n[s,t]\bigr|.
\]
Moreover, by the definition of $\tau^n$ \eqref{32} and $R^n$ \eqref{96},
we have $|R^n|\le c\varepsilon^n_M\vr^n<p_n<r_n$.
As a result \eqref{93} holds and the result follows.
\end{pf*}

For the proof of Lemma \ref{lem:aux3} we define $A_{kl}^n(t)$ to be the
number of class-$k$ customers entering pool $l$ by time $t$ and let its
centered and scaled version be given by
%
%
\begin{equation}
\hA_{kl}^n(t)=n^{-1/2}\bigl(A_{kl}^n(t)-
\mu_{kl}^n \xi _{kl}^*\nu_l
n^{1/2}t\bigr).\label{eq:pool_arrivals_scaled}
\end{equation}

\begin{pf*}{Proof of Lemma \protect\ref{lem:aux3}}
The proof of the lemma proceeds by induction on the class number.

\textit{Induction base}, $k=1$.
For this class all server pools $j\in\calJ(k)$ (if there are any) are
necessarily leafs of the tree.
Thus
if $j$ is such a pool, one has $\hB_{kj}^n=-\hI_j^n$ so that, by Lemma
\ref{lem:idleness}, $|\hB_{kj}^n|_{T^n}\leq c_1 p_n$.
Thus for $s<t\le\tau^n_{2,k}$,
%
\begin{eqnarray}\label{eq:leafpools}
\bigl|\mB_{kj}^n[s,t]\bigr| &\le& \bigl|\hB_{kj}^n\bigr|_{T^n}
\bigl(2+\mu_{kj}^n(t-s)\bigr)\leq c_2
\bigl(p_n+ n^{1/2}p_n(t-s)\bigr)
\nonumber
\\[-8pt]
\\[-8pt]
\nonumber
&\leq& c_2\bigl(p_n+ n^{1/2}s_n(t-s)
\bigr),
\end{eqnarray}
where we used the fact that $s_n/p_n\tinf$. Note that $c_2$ does not
depend on~$c^1_k$.

Consider next $j=\jbar(k)$. Using \eqref{31} we have
\begin{eqnarray*}
\bigl|\mB_{kj}^n[s,t]\bigr| &\leq&\bigl |\hQ_k^n[s,t]\bigr|+
\sum_{l\in\calJ(k)}\bigl|\mB_{kl}^n[s,t]\bigr|+\bigl\|
W^n[s,t]\bigr\|
\\
&\leq& c_3\bigl(p_n+ n^{1/2}s_n(t-s)+
\omega(t-s)\bigr).
\end{eqnarray*}
For $s<t\le\tau^n_{2,k}$,
we used Lemma \ref{lem:qtightnessuptostop} and \eqref{eq:leafpools}.
Note that
$c_3$ does depend on $c^1_k$.

\textit{Induction step}, $k>1$. Assume that the result of the
lemma holds for all $m<k$. Namely, if $m<k$ is such that $|\Del
_l^n|_{T^n}\leq c_l^1 p_n$ for all $l<m$, then \eqref{eq:mBincrements}
holds with~$m$ replacing $k$ and for all $0\leq s\leq t\leq T^n$ and
all $j\in\calJ(m)$ with a constant $\gamma_1$ that does not depend on
$c_m^1$. It holds for $j=\jbar(m)$ up to $\tau_{2,k}^n$ with a constant
$\gamma_2$ that does depend on $c_m^1$.

We will show that this holds for $k$. We thus assume that $|\Del
_m^n|_{T^n}\leq c_m^1$ for all $m<k$. By the induction assumption we
have the existence of a constant $c_1$, depending on $(c_m^1, m<k)$ but
not depending on $c_k^1$, such that \eqref{eq:mBincrements} holds for
all $m<k$ and all $l\sim m$. By the argument leading to \eqref
{eq:AthroughBbound} we have for all $s,t\leq T^n$ that
%
%
\begin{equation}
\label{eq:inductionA} \sum_{m,j<k}\bigl|
\hA_{mj}^n[s,t]\bigr|\leq c_4\bigl(p_n+
n^{1/2}s_n(t-s)+\omega(t-s)\bigr).
\end{equation}
Considering a pool $l\in\calJ(k)$, the idleness process satisfies
\begin{eqnarray*}
\hI_l^n(t)&=&\hI_l^n(s)-n^{-1/2}
\sum_{m\leq k}A_{ml}^n[s,t]
+n^{-1/2}\sum_{m\leq k}D_{ml}^n[s,t]
\\
&=&\hI_l^n(s)-\sum_{m\leq k}
\hA_{ml}^n[s,t]+ \sum_{m\leq k}
\mu_{ml}^n\int_s^t
\hB_{ml}^n(v)\,dv+\sum_{m\leq
k}V_{ml}^n[s,t].
\end{eqnarray*}
In turn,
%
%
\begin{eqnarray}
\bigl| \hA_{kl}^n[s,t]\bigr|\label{eq:interim111}&\leq&2\bigl|\hI_l^n\bigr|_{T^n}+
\sum_{m< k}\bigl|\hA _{ml}^n[s,t]\bigr|+
\sum_{m\leq k}\bigl|V_{ml}^n[s,t]\bigr|
\nonumber
\\[-8pt]
\\[-8pt]
\nonumber
&&{}+
\sum_{m\leq k}\biggl\llvert \mu_{ml}^n
\int_s^t\hB_{ml}^n(v)\,dv
\biggr\rrvert.
\end{eqnarray}
By the induction assumption \eqref{eq:mBincrements} holds for all
classes $m\in\I(l)$ so that, by Lemma~\ref{lem:idleness},
%
%
\begin{equation}
\bigl|\hI _l^n\bigr|_{T^n}\leq c_5p_n
\label{eq:interim100}
\end{equation}
for a constant that
does not depend on $c_k^1$. Also, $|\hB_{ml}^n|_{T^n}\leq s_n$ by
definition so that
%
%
\begin{equation}
\sum_{m\leq k}\biggl\llvert \mu_{ml}^n
\int_s^t\hB_{ml}^n(v)\,dv
\biggr\rrvert \leq c_7 n^{1/2}s_n(t-s).
\label{eq:interim110}
\end{equation}
Thus using \eqref{eq:inductionA}, \eqref{eq:interim100} and \eqref
{eq:interim110} in \eqref{eq:interim111} and applying the definition of
$\Om^n$ to bound the increments of $V^n$, we conclude that
\[
\bigl| \hA_{kl}^n[s,t]\bigr|\leq c_8
\bigl(p_n+ n^{1/2}s_n(t-s)+\omega(t-s)\bigr),
\]
where $c_8$ does not depend on $c_k^1$. By \eqref{eq:BkjConstruct} we
then have that
\begin{eqnarray*}
\bigl\llvert \mB_{kl}^n[s,t]\bigr\rrvert &\leq& \bigl|\hA
_{kl}^n[s,t]\bigr|+\bigl|V_{kj}^n[s,t]\bigr|
\\
&\leq&c_9\bigl(p_n+\omega(t-s)+ n^{1/2}s_n(t-s)
\bigr),
\end{eqnarray*}
where, as required, the constants do not depend on $c_k^1$ in the
definition of $\tau_{2,k}^n$ \eqref{eq:tau2kdefin}. This argument is
repeated for each $l\in\mL(k)$. The argument for the pool $\jbar(k)$
then follows exactly as in the induction basis and this concludes the
proof.
\end{pf*}

\begin{pf*}{Proof of Lemma \protect\ref{lem:qtightnessuptostop}}
By the definition of $\tau_{2,k}^n$, and since the jumps of both $\wh
Q^n$ and $\check X^n$ are $O(n^{-1/2})$ and recalling that $p_n \gg
n^{-1/2}$, we have that $|\Delta_k^n|_{\tau_{2,k}^n}\leq2c_k^1 p_n$.
The result is thus an immediate consequence
of Lemma~\ref{lem:omegatilde}.
\end{pf*}

\begin{pf*}{Proof of Lemma \protect\ref{lem:aux2}}
Recall that $k$ is fixed and $j=\jbar(k)$. Consider first the event
$\Om
^n_{k,U}$.
By \eqref{91}, $\Del^n_k$ remains positive on the interval $\calT
_n:=[\tau^n_{1,k},\tau^n_{2,k})$.
By part (ii) of the definition of the policy, during this interval
all service completions in pool $j$ are followed by admission to
service of customers from the
classes in $\mL(k)$. Also note that part (i) of this definition is
irrelevant during this time
interval because there are no idle servers at pool $j$ (indeed, the set
$\mathcal{K}_j$ is not empty on this interval; if there were any idle
servers in pool $j$
then they would be immediately assigned to customers of classes
in the set $\mathcal{K}_j$).
Thus, for $s,t\in\calT_n$, $s<t$, we have
\[
\sum_{l\in\mathcal{L}(k)}B_{lj}^n[s,t]= -
\sum_{l\in\mL(k)}D_{lj}^n[s,t]+\sum
_{l:l\sim j}D^n_{lj}[s,t].
\]
Using \eqref{10} and \eqref{21} we rewrite this as
\begin{eqnarray*}
\sum_{l\in\mathcal{L}(k)}B_{lj}^n[s,t]&=& -
\sum_{l\in\mathcal{L}(k)}\mu_{lj}^n\int
_s^t B_{lj}^n(v)\,dv-\sum
_{l\in\mL(k)} n^{1/2}V_{lj}^n[s,t]
\\
&&{} +\sum_{l:l\sim j}\mu_{lj}^n\int
_s^t B_{lj}^n(v)\,dv +\sum
_{l:l\sim j} n^{1/2}V_{lj}^n[s,t].
\end{eqnarray*}
Denote $\calL^c(k)=\{l\in\calI\dvtx l\sim j, l\notin\calL(k)\}$.
After centering and scaling we have
%
%
\begin{eqnarray}\label{eq:Bcentered}
\sum_{l\in\mathcal{L}(k)}\mB_{lj}^n[s,t]&=&
\sum_{l\in\calL^c(k)}\mu_{lj}^n
\xi_{lj}^*\nu_l(t-s) +\sum_{l:l\sim j}
\mu_{lj}^n\int_s^t
\hB_{lj}^n(v)\,dv
\nonumber
\\[-8pt]
\\[-8pt]
\nonumber
&&{}+\sum_{l\in\mL^c(k)}V_{lj}^n[s,t].
\end{eqnarray}
Note that $i:=\ibar(j)\in\mL^c(k)$ so that
the first term on the RHS of \eqref{eq:Bcentered} is bounded below by
%
%
\begin{equation}
\mu^n_{ij}\xi_{ij}^*\nu_j(t-s)\geq
c_1 n^{1/2}(t-s). \label{eq:interim66}
\end{equation}
Since $\|\hB^n\|_{T^n}\leq s_n$, the second term is bounded, in
absolute value, by\break
$c_2 n^{1/2}s_n(t-s)$. Since $s_n\to0$, this gives
\[
\sum_{l\in\mathcal{L}(k)}\mB_{lj}^n[s,t]
\ge c_3 n^{1/2}(t-s)-\bigl\|V^n[s,t]\bigr\|.
\]
Equation \eqref{eq:aux2_1} now follows by using the definition of $\Om
^n$ to bound the
increment of $V^n$.

Let us now consider $\Om^n_{k,D}$.
To prove \eqref{eq:aux2_2}, note by \eqref{91} that $\Del^n_k<0$ on the
time interval
$\calT^n$. Note that by the first part of the policy definition, new
\mbox{class-$k$} arrivals are not
sent to pool $j$ even if there are idle servers. Moreover, by the
second part,
upon each service completion, no class-$k$ customers are admitted into service
in pool $j$ during this time (since $\Del^n_k<0$).
Hence
\[
B_{kj}^n[s,t]=-D_{kj}^n[s,t]=-
\mu_{kj}^n\int_s^t
B_{kj}^n(v)\,dv + n^{1/2}V_{kj}^n[s,t],
\]
or, after scaling and centering,
\[
\mB_{kj}^n[s,t]=-\mu_{kj}^n
\xi_{kj}^*\nu_j (t-s) + V_{kj}^n[s,t].
\]
Exploiting again
the bounds for $V^n$ on $\Om^n$, we have
\eqref{eq:aux2_2}. This completes the proof.
\end{pf*}

\begin{pf*}{Proof of Lemma \protect\ref{lem:idleness}} Fix $j\in
\calJ$, $j\ne
j_0$. For $i=i_0$ we will use here, with some abuse of notation, $\J
(i_0)=\J(i_0)\setminus\{j_0\}$.

We start with an observation that relates the
condition of the lemma to the processes $\hA_{kj}^n$, $k\in\calI(j)$.
To that end, note that the process $\hB_{kj}^n$ satisfies the relation
%
%
\begin{eqnarray}
\label{eq:BkjConstruct} \hB_{kj}^n(t)&=&
\hB_{kj}^n(0)+n^{-1/2}\bigl(A_{kj}^n(t)-D_{kj}^n(t)
\bigr)
\nonumber
\\[-8pt]
\\[-8pt]
\nonumber
& =& \hB_{kj}^n(0)+\hA_{kj}^n(t)-
\mu_{kj}^n\int_s^t \hB
_{kj}^n(v)\,dv+V_{kj}^n(t).
\end{eqnarray}
Hence, using the definition of $\Om^n$ to bound the increment of $V^n$,
%
%
\begin{equation}
\bigl|\hA_{kj}^n[s,t]\bigr|\leq\bigl|\mB_{kj}^n[s,t]\bigr|+c_1
\bigl(p_n+\omega(t-s)\bigr). \label{eq:AthroughBbound}
\end{equation}
In view of this and the assumption of the lemma,
%
%
\begin{equation}
\label{98} \bigl|\hA_{kj}^n[s,t]\bigr|\leq c_2
\bigl(p_n+ n^{1/2}s_n(t-s)+\omega(t-s)\bigr), \qquad k
\in \calI(j).
\end{equation}
Let $i=\ibar(j)$ [note that $j\in\J(i)$] and define
\begin{eqnarray*}
\tau_{2}^n&:=&\inf \biggl\{s\geq0\dvtx\sum
_{l\in\J(i)}\hI _{l}^n(s)\geq \gamma
p_n \biggr\}\wedge T^n,
\\
\tau_{1}^n&:=&\sup \biggl\{s\leq\tau_2^n
\dvtx\sum_{l\in\J(i)}\hI _{l}^n(s)
\leq\gamma p_n/2 \biggr\}.
\end{eqnarray*}
Note that for all $s<t$, and each $l\in\J$,
%
%
\begin{eqnarray} \label{eq:hIdynamics}
\hI_{l}^n[s,t]&=&\frac{1}{ n^{1/2}}\sum
_{k:k\sim
l}\bigl(-A_{kl}^n[s,t]+D_{kl}^n[s,t]
\bigr)
\nonumber
\\[-8pt]
\\[-8pt]
\nonumber
&=&\sum_{k:k\sim l} \biggl(-\hA_{kl}^n[s,t]
-\mu_{kl}^n \int_{s}^t
\hB_{kl}^n(v)\,dv-V_{kl}^n[s,t]
\biggr).
\end{eqnarray}

On $[\tau_{1}^n,\tau_{2}^n)$ the tracking policy routes all class-$i$
arrivals to pools in the set $\J(i)$.
Hence on $[\tau_1^n,\tau_2^n)$, $\sum_{l\in\J
(i)}A_{il}^n[s,t]=A_{i}^n[s,t]$ and
%
%
\begin{eqnarray}
\label{eq:Aij_inc} \sum_{l\in\J(i)}
\hA_{il}^n[s,t]&=&n^{-1/2}\lambda_i^n(t-s)-
\sum_{l\in\J
(i)}\mu_{il}^n
\xi_{il}^* \nu_j(t-s)+\hA_i^n[s,t]
\nonumber
\\[-8pt]
\\[-8pt]
\nonumber
&\geq&\frac{c_3}{2} n^{1/2}(t-s)+\hA_i^n[s,t]
\end{eqnarray}
for a positive constant $c_3$.
The last inequality follows from the following observation: by \eqref
{08} and \eqref{15} we have that $\lambda_i^n=\lambda_i n+O( n^{1/2})=
\sum_{l\in\J(i)}\mu_{il}\xi_{kl}^*\nu_l n+\sum_{l\notin\J
(i)}\mu_{il}\xi
_{il}^*\nu_ln +O( n^{1/2})$.
If $i\neq i_0$, then by Assumption \ref{assn3}, $\xi_{i\jbar(i)}^*>0$.
If $i=i_0$, then $\xi_{ij_0}^*>0$. In either case, there exists
$l\notin\J(i)$ with $\xi_{il}^*>0$. Hence $n^{-1/2}(\lambda
_i^n-\sum_{l\in\J(i)}\mu_{il}\xi_{il}^*\nu_jn)\geq c_3 n^{1/2}$ for a positive
constant $c_3$ as required.

Using $\|\hB\|_{T^n}\leq s_n$, \eqref{eq:Aij_inc}, \eqref{98}
and the assumption of the lemma in \eqref{eq:hIdynamics}, we have
\begin{eqnarray*}
\sum_{l\in\J(i)}\hI_{l}^n[s,t]&
\leq& c_4\bigl( n^{1/2}s_n(t-s)-
n^{1/2}(t-s)\bigr) +\bigl|\hA_i^n[s,t]\bigr|+
\bigl\|V^n[s,t]\bigr\|
\\
&\le&-c_5 n^{1/2}(t-s) + c_6
\bigl(p_n+\omega(t-s)\bigr)
\end{eqnarray*}
for all $\tau_{1}^n\leq s\leq t\leq\tau_{2}^n$.
As in \eqref{eq:interim99} we have that $-c_5 n^{1/2}(t-s) +c_6\omega
(t-s)+c_6 g(n)\leq c_6 p_n$
for all $s\leq t$ and all sufficiently large $n$ and, in turn, that
$\sum_{l\in\J(i)}\hI_{l}^n[\tau_1^n,t]\leq c_{6}p_n$, for
$t\in[\tau_{1}^n,\tau_{2}^n)$. Note that $c_6$ does not depend on the
constant $\gamma$.
Hence $\gamma$ can be chosen in such a way that $\tau_{2}^n\geq T^n$.
Since $\hI^n\geq0$, $|\hI_j^n|_{T^n}\leq|\sum_{l\in\in\J(i)}\hI
_l^n|_{T^n}\leq c_6 p_n$. This completes the proof.
\end{pf*}

\begin{pf*}{Proof of Lemma \protect\ref{lem:XtoI}} Recall that
%
%
\begin{equation}
\theta'\hX^n(t)=\theta'
\hX^n(0)+\theta'W^n(t)+G^n(t)-R^n(t).
\label{eq:26}
\end{equation}

We treat separately each of the elements on the right-hand side above. First,
by~\eqref{eq:omegatilde2} we have on $\Om^n$ that $|\theta
'W^n(t)-\theta
' W^n(s)|\leq c_1(p_n + \omega(t-s))$. Recall that $\alpha_\omega
>2\alpha_g$
so that $\omega(\vartheta_n)< n^{-(1/2\alpha_\omega-\delta\alpha
_\omega)}\leq
n^{-\alpha_g+\delta}=p_n n^{\delta}=q_n$. Thus
%
%
\begin{equation}
\bigl|\theta'W^n(t)-\theta'W^n(s)\bigr|
\leq c_2 q_n.\label{eq:24}
\end{equation}
Next, by
the definition of $T^n$, using \eqref{32}, we have that
%
%
\begin{equation}
\bigl|R^n\bigr|_{T^n}\leq J\|\theta\|\varepsilon_M^n
\varrho^n\leq p_n\label {eq:25}.
\end{equation}
Recalling that $\theta_i\mu_{ij}=z_j^*$ for all $i\sim j$
and that $\sum_{i}\tB_{ij}^n=-I_j^n$ for all $j\in\J$ we have that
%
%
\begin{equation}
G^n(t)=-\sum_{i,j}\theta_i
\mu_{ij}\int_0^t\tB
_{ij}^n(u)\,du=\sum_{j}z_j^*
\int_0^t I_j^n(u)\,du,
\label{eq:Gidentity}
\end{equation}
from which we also
see that $G^n$ is nondecreasing. By Proposition \ref
{prop:localadjustment1},\break  $|\Del_k^n|_{T^n}\leq c_k^1p_n$ for all $k\in
\I
_{-i_0}$ so that, by Lemma \ref{lem:idleness}, $|\hI_j^n|_{T^n}\leq c_3
p_n$ for all $j\neq j_0$. In turn, $\int_s^t I_j^n(u)\,du\leq c_4 p_n
n^{1/2}(t-s)$ for all $j\neq j_0$ and for all $s,t\leq T^n$,
\[
\biggl|G^n(t)-G^n(s)-z_{j_0}^*\int
_s^t I_{j_0}^n(u)\,du \biggr|\leq
\sum_{j\neq j_0}z_j^*\int_s^t
I_j^n(u)\,du\leq c_5 n^{1/2}p_n
(t-s).
\]
Letting $t-s=\vartheta^n$ proves the first part of the lemma. The
second part then follows immediately from \eqref{eq:24}, \eqref{eq:25}
and the fact that $G^n$ is nondecreasing.
\end{pf*}

\begin{pf*}{Proof of Lemma \protect\ref{lem:cXoscillation}}
Throughout we fix
$s,t\leq T^n$ as in the statement of the lemma. We argue separately for
two cases according to whether there exists $u\in[s,t)$ such that
$\theta'\hX^n(u)\geq\bar\kappa_n$.

Suppose first that $\theta'\hX^n(u)<\bar\kappa_n$ for all $u\in
[s,t)$. By the properties of the functions $f^n$ (see Remark \ref
{rem:qnproperties}) and since $\kappa_n\ll q_n$ we have here that $\cX
_i^n[s,t]\leq(\bI\theta_i)^{-1}\kappa_n\ll q_n$ as required.

To treat the other case we establish first the following claim. Fix
$\beta>0$, then for all sufficiently large $n$, if $u\leq T^n$ has
$\theta'\hX^n(u)\geq\beta\bar\kappa_n$, then $\hQ_i^n(u)>0,
i\in\I$.

To see this, note by Proposition \ref{prop:localadjustment1} and the
properties of $f^n$ (see Remark~\ref{rem:qnproperties}) that $\theta
_i\hQ_i^n(u)\geq\theta_if_i^n(\theta'\hX^n(u))-c_1 p_n\geq(\bI
\theta
_i)^{-1}\kappa_n-c_1p_n>0$ for all $i\in\I_{-i_0}$ where we use the
fact that $p_n\ll\kappa_n$. For $i=i_0$,
\begin{eqnarray*}
\theta_{i_0}\hQ_{i_0}^n(u)&=&
\theta'\hX^n(u)-\sum_{i\neq i_0}
\theta_i\hQ_i^n(u)-\sum
_{i,j}\theta_i\hB_{ij}^n(u)
\\
&\geq& \theta'\hX^n(u)-\sum
_{i\neq i_0}\theta_i\cX_i^n(u)-c_2
s_n-c_3p_n
\\
&=& \theta_{i_0}\cX_{i_0}^n(u)-c_4
s_n>0.
\end{eqnarray*}
The first inequality follows from Proposition \ref
{prop:localadjustment1} and the fact that $\|\hB^n\|_{T^n}\leq s_n$ by
the definition of $T^n$. The second inequality follows from $p_n\ll
s_n$ and from the fact that $\theta'\hX^n=\theta'\cX^n$ whenever
$\theta
'\hX^n\geq0$. The last inequality then follows from the definition of
$\cX^n$, the fact that $f_{i_0}^n(x)\geq(\bI\theta_{i_0})^{-1}(\bar
\kappa_n\wedge x)$ for all $x\geq0$ (see Remark~\ref
{rem:qnproperties}) and recalling that $s_n \ll\bar\kappa_n$.

Having the above we proceed to consider the case in which $\theta'\hX
^n(u)\geq\bar\kappa_n$ for some $u\in[s,t)$. Let
\[
\tau_1^n=\sup\bigl\{\eta\leq u\dvtx
\theta'\hX^n(u)\leq\tfrac{1}{2}\bar
\kappa_n\bigr\}
\]
and
\[
\tau_2^n=\inf\bigl\{\eta\geq u\dvtx
\theta'\hX^n(u)\leq\tfrac{1}{2}\bar
\kappa_n\bigr\} \wedge t,
\]
where we set $\tau_1^n=s$ if $\theta'\hX^n(\eta)\geq\bar\kappa
_n$ for
all $\eta\in[s,u)$. Since the jumps of $\theta'\hX^n$ are of size
$O(n^{-1/2})$ we must have that, if $\tau_2^n<t$ or $\tau_1^n>s$, then
$\tau_1^n<u<\tau_2^n$. Setting $\beta=1/2$ in the argument above we
have that $\hQ_i^n>0$ on $[\tau_0^n,\tau_1^n)$ so that, by Remark
\ref
{rem:workconservation} $\hI_j^n=0, j\in\J$ on this interval. By
equation \eqref{eq:Gidentity} we then have that $G^n(\tau
_2^n)-G^n(\tau
_1^n)=0$. Using \eqref{eq:26} and \eqref{eq:25} we then have that
%
%
\begin{eqnarray}
\bigl|\theta'\hX^n\bigl(\tau_1^n
\bigr)-\theta'\hX^n\bigl(\tau_2^n
\bigr)\bigr|&\leq&\bigl|\theta 'W^n\bigl(\tau _2^n
\bigr)-\theta W^n\bigl(\tau_1^n
\bigr)\bigr|+c_4 p_n
\nonumber
\\[-8pt]
\\[-8pt]
\nonumber
&\leq & c_5
\bigl(p_n+\omega\bigl(\tau _2^n-\tau
_1^n\bigr)\bigr).
\end{eqnarray}
It then follows, as in the beginning of the proof of Lemma
\ref{lem:XtoI}, that
%
%
\begin{equation}
\bigl|\theta'\hX^n\bigl[\tau_1^n,
\tau_2^n\bigr]\bigr|\leq c_6 q_n.
\label{eq:31}
\end{equation}
Since $q_n\ll\bar\kappa_n$ we conclude that
$\theta'\hX^n\geq\frac{3}{4}\bar\kappa_n$ for all $u\in[\tau
_1,\tau
_2^n)$ which contradicts $s<\tau_1^n<\tau_2^n<t$. We conclude that
$\tau
_1^n=s$ and $\tau_2^n=t$ and, by \eqref{eq:30}, that $\theta'\hX
^n[s,t]\leq c_7q_n$. From the local Lipschitz continuity of $f^n$ it
then follows that $\cX_i^n[s,t]\leq c_8q_n$ as required.
\end{pf*}

\begin{pf*}{Proof of Lemma \protect\ref{lem:Bbound}}
We use the notation $\bar{w}_t:=\bar{w}_{[0,t]}$ with the latter
defined in
\eqref{eq:wdefin}. Recall also that the index set $\calI$ is identified
with $\{1,\ldots,\bI\}$ and define
%
%
\begin{equation}
b_k(n)=n^{k\delta/\bI}p_n,\qquad k\in\{0,1,\ldots,\bI\}.
\label{eq:bkdefin}
\end{equation}
Note that $b_{\bI}(n)=n^{\delta}p_n=q_n$.

Fix $k\in\I_{-i_0}$ and let $j=\jbar(k)$. By \eqref{31} and using
$\hat
X^n(0)=0$ we can write
%
%
\begin{equation}
\hB_{kj}^n(t)=-\mu_{kj}^n\int
_0^t \hB_{kj}^n(v)\,dv+F^n(t),\qquad
t\ge0, \label{eq:Bkjrep}
\end{equation}
where $F^n=F_1^n+F_2^n$,
\[
F_1^n(t)=-\wh Q^n_k(t)+W_k^n(t),\qquad
F_2^n(t)=-\sum_{l\in\calJ(k)}
\mu_{kl}^n\int_0^t\wh
B^n_{kl}(s)\,ds -\sum_{l\in\calJ(k)}
\hB_{kl}^n(t).
\]
The proof is based on the following estimate; see Lemma 3.4 of \cite{atasol}
and its proof. Let $X$ be the unique solution to the integral equation
\[
X(t)=-\mu\int_0^tX(s)\,ds+F(t),\qquad t\ge0,
\]
with $\mu>0$ and data $F\dvtx[0,\iy)\to\D$.
Then given $u>0$ and $\vartheta\in(0,u)$,
%
%
\begin{equation}
\label{97} |X|_u\le2|F|_ue^{-\mu\vartheta}+\bar
w_u(F,\vartheta).
\end{equation}
Thus in view of \eqref{eq:Bkjrep}, one can bound $\hat B^n_{kj}$ by
suitably estimating $F^n$.

If $\sum_{l\in\calJ(k)}|\hB_{kl}^n|_{T^n}\leq\beta b_{k-1}(n)$ for
some $\beta$, then given $\vartheta>0$,
%
%
\begin{equation}
\bar{w}_{T^n}\bigl(F_2^n,\vartheta\bigr)\leq
c_1b_{k-1}(n) \bigl(1+\vartheta n^{1/2}\bigr),
\label{eq:F2modulus}
\end{equation}
where $c_1$ does not depend on $n$ or $\vartheta$. By Proposition \ref
{prop:localadjustment1} we have that $|\hQ_k^n-\cX_k^n|_{T^n}\leq c_2p_n$.
Combined with Lemma \ref{lem:cXoscillation}, and letting $\vartheta
=\vartheta_n$ this gives
\[
\bar{w}_{T^n}\bigl(F_1^n,\vartheta_n
\bigr)\leq c_3 p_n
\]
for some constant $c_3$ (not depending on $n$). Noting that $
n^{1/2}\vartheta_n=n^{\delta}$, we get
%
%
\begin{equation}\label{eq:Fmodul}
\bar{w}_{T^n}\bigl(F^n,\vartheta_n\bigr)\leq
 c_4 b_{k-1}(n)n^{\delta}=c_4
b_k(n),
\end{equation}
where we used the fact that $ n^{1/2}b_{k-1}(n)\vartheta_n=n^{\delta
}b_{k-1}(n)=b_k(n)$ [recall
that $b_k(n)=n^{({|k|}/{\bI})\delta}p_n$]. Noting $F^n(0)=0$ by our
assumptions and using \eqref{eq:omegatilde3} and~\eqref{eq:omegatilde2},
we also have
%
%
\begin{equation}
\bigl|F^n\bigr|_{T^n}\leq c_5\bigl(r_n+1+b_{k-1}(n)
n^{1/2}\bigr). \label{eq:Fbound}
\end{equation}
Thus
%
%
\begin{equation}
\bar{w}_{T^n}\bigl(F^n,\vartheta_n\bigr)\leq
c_6b_k(n),\label{eq:Fmodul2}
\end{equation}
and, in turn,
\[
\bigl|\hB_{kj}^n\bigr|_{T}\leq\bigl|F^n\bigr|_{T^n}e^{-\mu_{kj}^n
\vartheta_n}+
\bar{w}_t\bigl(F^n,\vartheta_n\bigr)
\bigl(1-e^{-\mu_{kj}^n\vartheta_n}\bigr).
\]
Since $ n^{1/2}\vartheta_n=n^{\delta}$ we have that $\mu
_{kj}^n\vartheta
_n\geq c_7n^{\delta}$. Further, since
$b_{k-1}(n)\rightarrow0$, we have that $ n^{1/2}b_{k-1}(n)\leq c_8
n^{1/2}$ so that
$ n^{1/2}z(n)e^{-n^{\delta}}\rightarrow0$. We conclude that
%
%
\begin{equation}
\bigl\|\hB_{kj}^n\bigr\|_{T^n}\leq2\bar{w}_{T^n}
\bigl(F^n,\vartheta_n\bigr)\leq c_9b_k(n),
\label{eq:Bboundinducstep}
\end{equation}
provided that
%
%
\begin{equation}
\sum_{l\in\calJ(k)}\bigl|\hB _{kl}^n\bigr|_{T^n}
\leq\beta b_{k-1}(n)\label{eq:Bboundinducasum}
\end{equation}
for
some $\beta$.

The requirement \eqref{eq:Bboundinducasum} holds trivially if $k$ is a
leaf of the tree, in which case the set $\J(k)$ is empty. It also holds
if all pools in $\J(k)$ are leafs of the tree by Proposition \ref
{prop:localadjustment1} and Lemma \ref{lem:idleness} because in that
case $|\hB_{kj}|_{T^n}=|\hI_j^n|_{T^n}$ for all $j\in\J(k)$. Thus the
fact that it holds for all $k\neq i_0$ and $j=\jbar(k)$ now follows by
induction on the class number using \eqref{eq:Bboundinducasum} and
\eqref{eq:Bboundinducstep}.

To bound $\hB_{kj}^n$ for $j\neq\jbar(k)$, by identity \eqref{12} we
have that
$|\hB_{kj}^n|_{T^n}\leq|\hI_j^n|_{T^n}+\sum_{l\in I(j)}|\hB
_{lj}|_{T^n}$. The first part of this proof guarantees that\break  $\sum_{l\in
I(j)}|\hB_{lj}|_{T^n}\leq c_{10}q_n$. From Proposition \ref
{prop:localadjustment1} and Lemma \ref{lem:idleness} it follows that
$|\hI_j^n|_{T^n}\leq c_{10}p_n$. This completes the proof.
\end{pf*}

\begin{pf*}{Proof of Lemma \protect\ref{lem:oscillationI}}
Let $s,t\leq T^n$ be an interval as in the statement of the lemma. In
particular, $\hI_{j_0}^n>0$ on $[s,t)$ so that, by Remark \ref
{rem:workconservation}, $\hQ_{i_0}^n=0$ on $[s,t)$ and $\theta
_{i_0}\hX
_{i_0}^n[s,t]=\sum_{j}\hB_{i_0j}^n[s,t]$. Thus
\[
\theta_{i_0}\hB_{i_0j_0}^n[s,t]=
\theta_{i_0}X_{i_0}^n[s,t]-\theta _{i_0}
\sum_{j\neq j_0}\hB_{i_0j}^n[s,t]\geq
\theta_{i_0}\hX _{i_0}^n[s,t]-c_1
q_n,
\]
where the inequality follows from Lemma \ref{lem:Bbound}.

By Lemma \ref{lem2} and the identity $\hI_{j_0}^n=-\sum_{i}\hB
_{ij_0}^n$ it further holds that
\[
\bigl|\hI_{j_0}[s,t]+\hB_{i_0j_0}^n[s,t]\bigr|\leq\sum
_{i\neq i_0}\bigl|\hB _{ij_0}^n\bigr|_{T^n}
\leq c_2 q_n,
\]
so that
%
%
\begin{equation}
\theta_{i_0}\hI_{j_0}[s,t]\leq-\theta_{i_0}X_{i_0}^n[s,t]+c_4
q_n.\label{eq:23}
\end{equation}
The lemma would then follow from Lemma \ref{lem:XtoI} provided that
%
%
\begin{equation}
\theta_{i_0}\hX_{i_0}^n[s,t]\geq
\theta'\hX^n[s,t]- c_{5} q_n.
\label {eq:30}
\end{equation}
We next prove \eqref{eq:30}. By Lemma \ref{lem:Bbound}, $\|\hQ
_{-i_0}^n-\hX_{i_0}^n\|_{T^n}=\sum_{i\in\I_{-i_0},j\in\J}|\hB
_{ij}^n|_{T^n}\leq c_6 q_n$ and by Proposition
\ref{prop:localadjustment1},
$\|\hQ_{-i_0}^n-\cX_{-i_0}^n\|_{T^n}\leq c_7 p_n$. In turn,
\[
\bigl\|\hX_{-i_0}^n-\cX_{-i_0}^n
\bigr\|_{T^n}\leq\bigl\|\hQ_{-i_0}^n-\hX _{-i_0}^n
\bigr\| _{T^n}+\bigl\|\hQ_{-i_0}^n-\cX_{-i_0}^n
\bigr\|_{T^n}\leq c_8 q_n.
\]
Using the identity $\theta_{i_0}\hX_{i_0}^n[s,t]=\theta'\hX
^n[s,t]-\sum_{i\neq i_0}\hX_i^n[s,t]$ we then have that
%
%
\begin{equation}
\biggl|\theta_{i_0}\hX _{i_0}^n[s,t]-
\theta'\hX^n[s,t]+\sum_{i\neq i_0}
\cX_i^n[s,t]\biggr|\leq c_9 q_n.
\label{eq:27}
\end{equation}
By Lemma \ref{lem:cXoscillation}, $\sum_{i\neq
i_0}|\cX_i^n[s,t]|\leq c_{10}q_n$ for all $s,t\leq T^n$ with
$|t-s|\leq
\vartheta_n$ which proves \eqref{eq:30} and thus completes the proof of
the lemma.
\end{pf*}

%



\printaddresses

\end{document}